\documentclass{article}

\usepackage[english]{babel}

\usepackage[letterpaper,top=2cm,bottom=2cm,left=3cm,right=3cm,marginparwidth=1.75cm]{geometry}

\usepackage{amsmath}
\usepackage{graphicx}
\usepackage[colorlinks=true, allcolors=blue]{hyperref}
\usepackage{subfigure}
\usepackage{wasysym}
\usepackage{caption}
\usepackage{lipsum}
\usepackage{booktabs}
\usepackage{amssymb}
\date{}
\begin{document}

\title{Periodic micromagnetic finite element method}

\author{Fangzhou Ai \and Jiawei Duan\thanks{Co-author} \and Vitaliy Lomakin\thanks{Corresponding Author. Email: vlomakin@ucsd.edu, also at Program in Materials Science and Engineering, University of California, San Diego, La Jolla 92093, U.S.A}}

\maketitle

\begin{center}
    Department of Electrical and Computer Engineering, University of California San Diego, \\
    La Jolla 92093, U.S.A \\
    Center for Memory and Recording Research, University of California San Diego, \\
    La Jolla 92093, U.S.A \\
\end{center}

\begin{abstract}

Periodic micromagnetic finite element method (PM-FEM) is introduced to solve periodic unit cell problems using the Landau-Lifshitz-Gilbert equation. PM-FEM is applicable to general problems with 1D, 2D, and 3D periodicities. PM-FEM is based on a non-periodic FEM-based micromagnetic solver and extends it in several aspects to account for periodicities, including the computation of exchange and magnetostatic fields. For the exchange field, PM-FEM modifies the sparse matrix construction for computing the Laplace operator to include additional elements arising due to the periodicities. For the magnetostatic field, the periodic extensions include modifications in the local operators, such as gradient, divergence, and surface magnetic charges as well as the long-range superposition operator for computing the periodic scalar potential. The local operators are extended to account for the periodicities similar to handling the Laplace operator. For the long-range superposition operator, PM-FEM utilizes a periodic Green's function (PGF) and fast spatial convolutions. The PGF is computed rapidly via exponentially rapidly convergent sums. The spatial convolutions are accomplished via a modified fast Fourier transform based adaptive integral method that allows calculating spatial convolutions with non-uniform meshes in $O(N\log N)$ numerical operations. PM-FEM is implemented on CPU and GPU based computer architectures. PM-FEM allows efficiently handling cases of structures contained withing the periodic unit cell touching or not touching its boundaries as well as structures that protrude beyond the unit cell boundaries. PM-FEM is demonstrated to have about the same or even higher performance than its parent non-periodic code. The demonstrated numerical examples show the efficiency of PM-FEM for highly complex structures with 1D, 2D, and 3D periodicities.

\end{abstract}

\section{Introduction} \label{sec:intro}
Periodic structures are common in many fields of study, such as crystals\cite{vasp,QE-2009}, molecular dynamics\cite{HUNENBERGER199969, Weber2000}, electromagnetics\cite{6955453, mailloux2017phased, kalkstein1971green, sirenko2010modern, peng1975theory} among others. Periodic structures are also of a high importance in Micromagnetics \cite{osti_1349900}. For magnetic structures, periodicity can be used to mimic infinite domains, e.g., infinite wires, films, or bulk. Periodicity can be used to account for spin wave propagation. It also can be used to approximate structures that are large in one of their dimensions. 

There is a range of methods used in various fields of study to account for periodicity \cite{1504955, 774139, https://doi.org/10.1029/2004RS003171, 7778376, r1_1, r1_2, r1_3, r1_4, AI2024109291}, or specific algorithms customized for Micromagnetics\cite{WANG201084, fdtd, WYSOCKI2017274, PhysRevB.87.174422,  PhysRevB.87.174422}. In Micromagnetics, approaches accounting for periodicity are based on finite difference method (FDM) based domain discretization schemes with uniform grids\cite{fdtd}. For a uniform grid with $N$ grid points, the magnetostatic field accounting for periodicity can accomplished using the Fast Fourier Transform (FFT) method in $O(N\log N)$ computational cost. Even for the FDM, one needs to be careful in the approaches of computing the periodic superposition kernel, referred to as periodic Green's function (PGF) as such computational can become slow or even lead to improper results\cite{PhysRevB.57.14332}. The situation is more complicated when using micromagnetic codes based on the finite element method (FEM). For FEM, the structure to be modelled is meshed into a generally non-uniform mesh, which is often tetrahedral or hexahedral. Such meshing allows for a great flexibility in modeling complex materials and devices. However, the mesh non-uniformity makes it more complicated to account for general periodicities. Currently, there are no reported numerical methods or codes for micromagnetics that can handle general periodic structures.

Here, we introduce a formulation, referred to as periodic micromagnetic finite element method (PM-FEM) to handle periodicities in complex micromagnetic simulations using FEM. PM-FEM offers a unified framework incorporating 1D, 2D, and 3D periodicities. For the computation of the exchange field, the PM-FEM modifies the exchange matrix to make use of the corresponding (1D, 2D, 3D) PBC. For the computation of the magnetostatic field, the PM-FEM adapts the Green's function to PGF with a modified precorrected Fast Fourier Transform or adaptive integral method to calculate the periodic scalar potential (PSP) followed by using a periodized gradient operator\cite{huebner2001finite}. Additionally, PM-FEM allows accounting for periodic unit cells of different configurations, including cases of touching and non-touching boundaries as well as cases of contained and protruding unit cells. The paper is organized as follows. Section~\ref{sec:formula} presents the formulation of the problem. Section~\ref{sec:imple} presents the implementation details for the localized and non-localized interactions. Section~\ref{sec:res} shows results obtained using the PM-FEM method, including its performance and error analysis. Finally, Sec.~\ref{sec:con} presents summary and conclusions.

\section{Formulation} \label{sec:formula}

\subsection{Continuous formulation} \label{sec:formula:cont}

\begin{figure}[htbp] \label{fig:categ}
    \centering
    \subfigure[]{
    \includegraphics[height=0.11\linewidth]{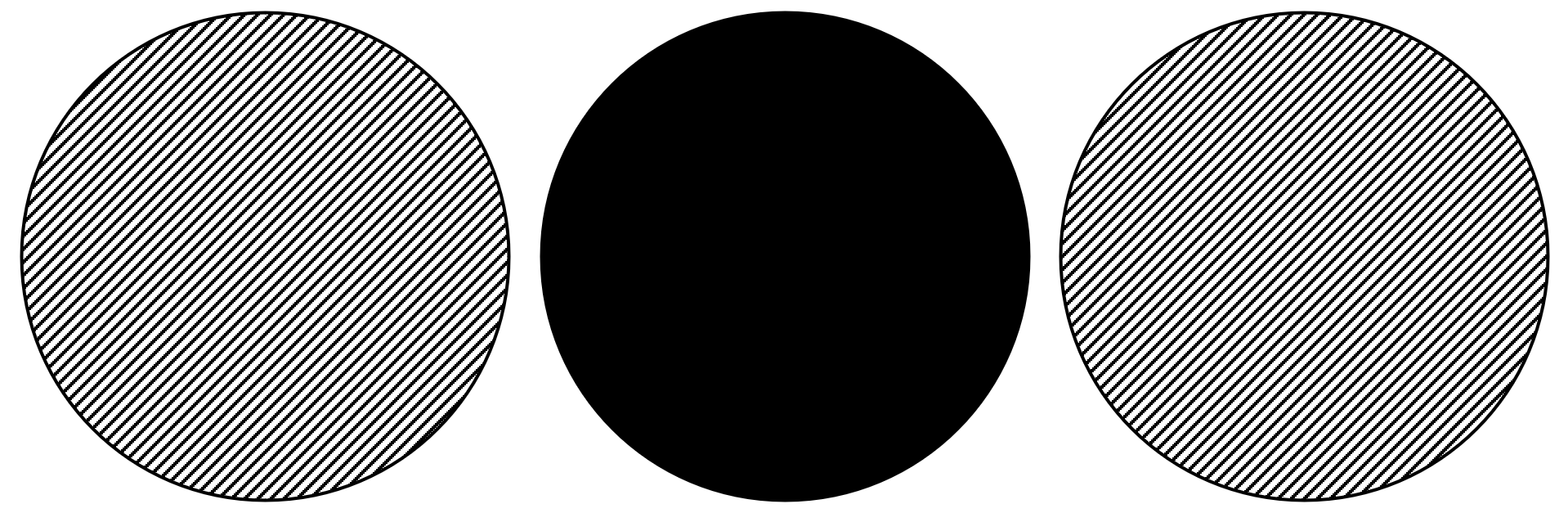} \label{fig:nt-np}
    }
    \quad
    \subfigure[]{
    \includegraphics[height=0.1\textwidth]{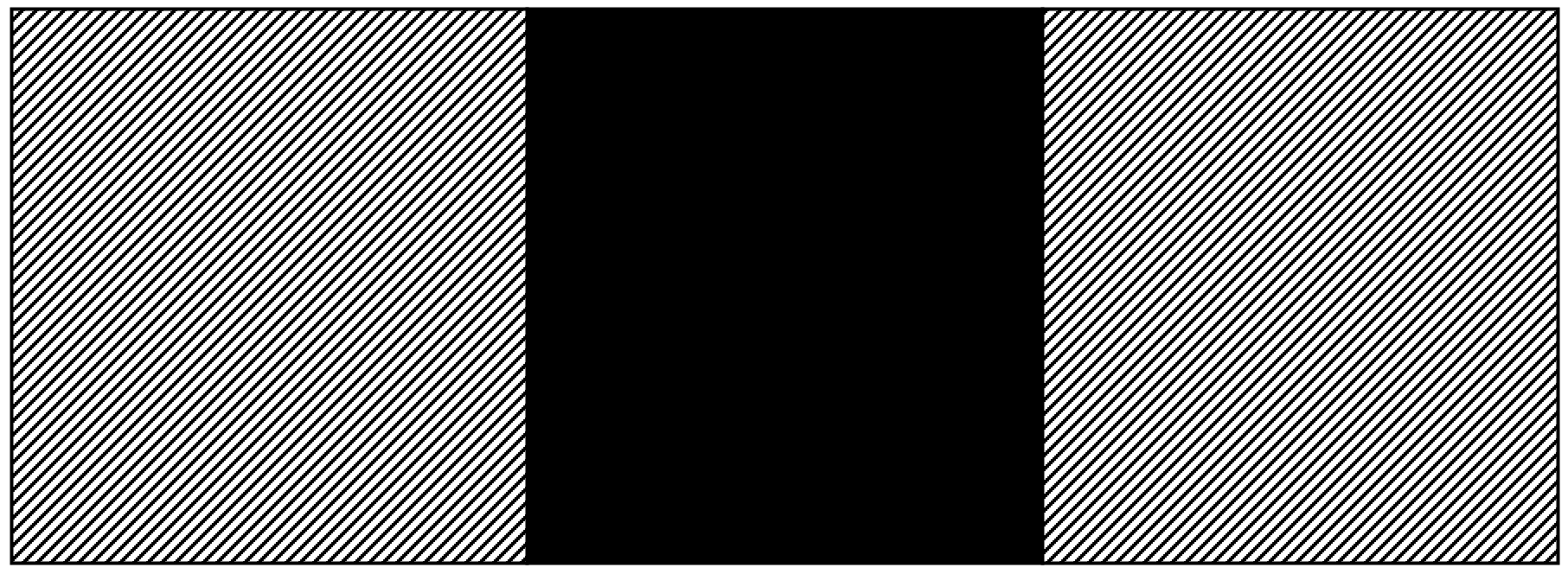} \label{fig:t-np}
    }
    \newline
    \subfigure[]{
    \includegraphics[height=0.1\textwidth]{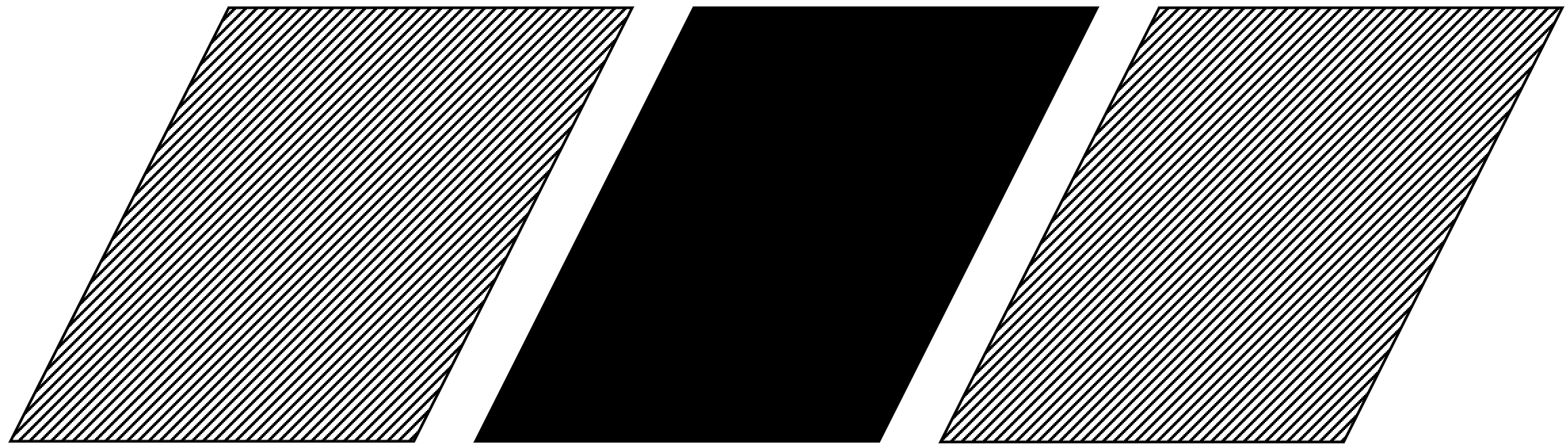} \label{fig:nt-p}
    }
    \quad
    \subfigure[]{
    \includegraphics[height=0.1\textwidth]{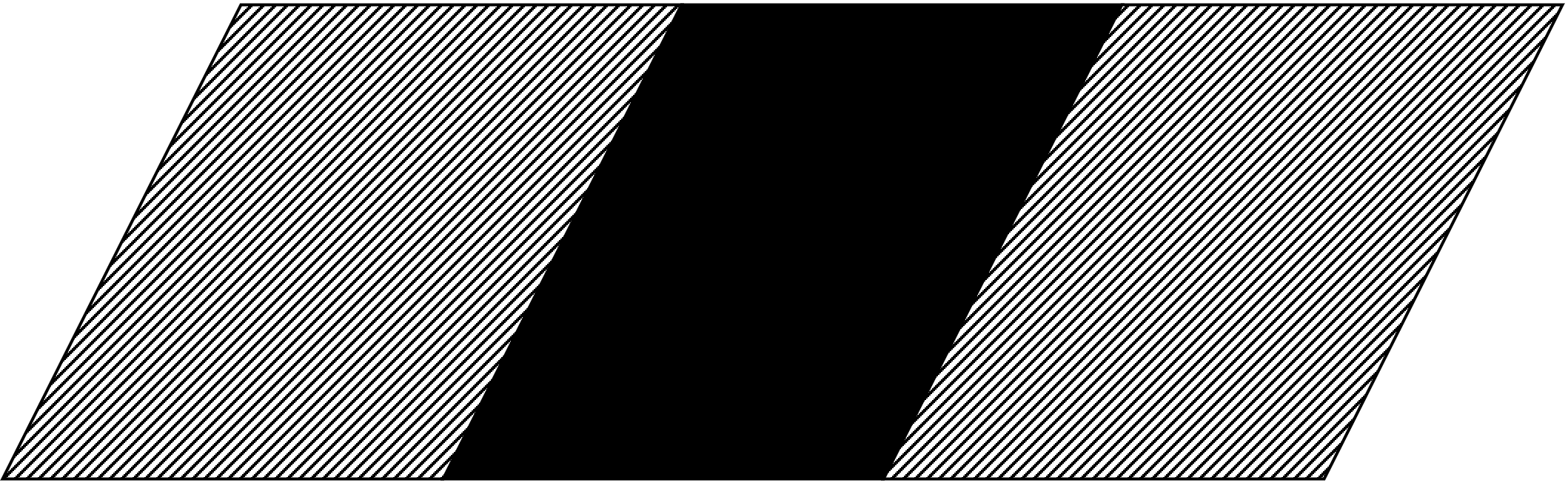} \label{fig:t-p}
    }
    \caption{Categories of 1D PBC unit cell (solid black) and its nearest images (shadowed), (a) non-touching, non-protruding case, (b) touching, non-protruding case, (c) non-touching, protruding case and (d) touching, protruding case.}
    \label{fig:categ}
\end{figure}

We consider a magnetic structure characterized by the magnetization $\mathbf{M}(\mathbf{r})$, where $\mathbf{r}$ is the location vector defined in a periodic unit cell with 1D, 2D, or 3D PBC (Fig. 1). For the 3D PBC the periodic unit cell has periodic dimensions (periods) of $L_x$, $L_y$, and $L_z$ (Fig. 1a). For the 2D PBC, the periodic dimensions are $L_x$ and  $L_y$, whereas the other ($z$) dimension is non periodic. For the 1D PBC, the periodic dimension is $L_x$, whereas the other ($y$ and $z$) dimensions are non-periodic. The magnetization satisfies the following periodic condition: 
\begin{subequations} \label{eq:pm_cont}
    \begin{align} 
    &\mathbf{M}(\mathbf{r}+L_x\hat{x})=\mathbf{M}(\mathbf{r}),\\ 
    &\mathbf{M}(\mathbf{r}+L_x\hat{x}+L_y\hat{y})=\mathbf{M}(\mathbf{r}),\\ 
    &\mathbf{M}(\mathbf{r}+L_x\hat{x}+L_y\hat{y}+L_z\hat{z})=\mathbf{M}(\mathbf{r})
\end{align}
\end{subequations}
where $L_x$, $L_y$, and $L_z$ are the periods for the 1D, 2D, and 3D PBC, respectively. There are two orthogonal PBC categories that require a special attention in numerical treatment, touching-(T-)/non-touching-(NT-)PBC and protruding (P-)/non-protruding(NP-)PBC. For T-PBC, such cases occur when for a point at location $\mathbf{r}_i,i\in[1,N]$, one of its periodic images $\mathbf{r}_i+\gamma L_x\hat{x}+\zeta L_y\hat{y}+\xi L_z\hat{z}$ can still be found within the $0^{th}$ unit cell, namely we can find another point located at $\mathbf{r}_j = (x_i+\gamma L_x, y_i+\zeta L_y, z_i+\xi L_z), j\in[1,N], j\neq i$. Here, $\gamma$, $\zeta$ and $\xi$ are zero if no periodicity is present at the $\hat{x}$, $\hat{y}$ or $\hat{z}$ direction, and they can be $1$ if a periodicity is present. We also have $\gamma + \zeta + \xi > 0$ to avoid $j=i$, if for any $\mathbf{r}_i$, we can't find the corresponding $\mathbf{r}_j$, then it is NT-PBC. Regarding of the P-PBC, it is defined as the unit cell is tightly confined by the periodic length, namely the unit cell size is strictly less or equal than periodic length on all directions, and NP-PBC is defined if the condition does not hold. These two categories are orthogonal to each other, and all periodic problems can be classified as arbitrary combination of these two categories, e.g. T-NP-PBC, NT-P-PBC etc. Here we show some 1D PBC examples for illustration purpose. Fig.~\ref{fig:nt-np} shows a NT-NP-PBC unit cell, it can be represented as an infinite array of spheres, Fig.~\ref{fig:t-np} shows a T-NP-PBC, it can be used to model an infinite wire, and Fig.~\ref{fig:nt-p} and Fig.~\ref{fig:t-p} shows the NT-P-PBC and T-P-PBC respectively. All these cases can be met in micromagnetic simulations and need to be addressed. Specifically, T-PBC (either protruding or non-protruding) requires a special treatment for handling relevant differential operators whereas the protruding PBC (either touching or non-touching) requires special handling of the long-range superposition integrals / sums.

The magnetization dynamics is governed by the Landau-Lifshits-Gilbert equation:
\begin{equation} \label{eq:llg}
\frac{\partial \mathbf{M}}{\partial t} = -\frac{\gamma_{\text{gr}}}{1+\alpha^2}\left(\mathbf{M}\times \mathbf{H}_\text{eff} + \frac{\alpha}{M_s}\mathbf{M}\times\mathbf{M}\times\mathbf{H}_\text{eff}\right),
\end{equation}
where $\gamma_{\text{gr}}$ is the gyromagnetic ratio, $\alpha$ is the damping constant, $M_s$ is the saturation magnetization, and 
\begin{equation} \label{eq:heff}
\mathbf{H}_\text{eff} = \mathbf{H}_\text{ex} + \mathbf{H}_\text{ms} + \mathbf{H}_\text{an} + \mathbf{H}_\text{ap}
\end{equation}
is the effective magnetic field that includes the exchange, magnetostatic, anisotropy, and applied field components. The effective field also can include additional components, such as spin transfer torque \cite{RALPH20081190}, spin orbit torque \cite{shao2021roadmap}, eddy currents \cite{hrkac2005three}, magnetostriction \cite{ANDRA20018892}, etc. The effective field and all its components satisfy the periodicity conditions, as in Eq. \eqref{eq:pm_cont}.

\subsection{Discrete formulation} \label{sec:formula:discr}

For numerical calculations, we utilize the finite element method (FEM), in which the unit cell is meshed into a set of elements. We choose tetrahedrons as the discretization elements for a broad availability of meshers but also hexahedrons can be a good choice. The magnetization is represented as  
\begin{equation} \label{eq:M_phi}
\mathbf{M} = \sum_{n=1}^N \mathbf{M}_n\varphi_n(\mathbf{r}),
\end{equation}
where $\mathbf{M}_n = \mathbf{M}(\mathbf{r}_n)$ and $\varphi_n(\mathbf{r})$ are basis functions chosen as linear interpolatory polynomials. The solutions for the magnetization $\mathbf{M}_n$ are obtained at points $\mathbf{r}_n$ and the magnetization at any other point can be found via Eq.~\eqref{eq:M_phi}.

The discretization in Eq.~\eqref{eq:heff} leads to a a unit cell that contains $N$ coinciding source (magnetizations) / observer (effective fields) points. Each point is located at $\mathbf{r}_n=(x_n, y_n, z_n), n\in [1,N]$ and has magnetization $\mathbf{M}(\mathbf{r}_n)=(M_n^x, M_n^y, M_n^z)$ with an infinite number of images in the corresponding periodic dimensions. The magnetization in each image satisfies the following periodic conditions:
\begin{subequations} \label{eq:pm}
    \begin{align} 
    &\mathbf{M}(\mathbf{r}_n+i_xL_x\hat{x})=\mathbf{M}(\mathbf{r}_n),\\ 
    &\mathbf{M}(\mathbf{r}_n+i_xL_x\hat{x}+i_yL_y\hat{y})=\mathbf{M}(\mathbf{r}_n),\\ 
    &\mathbf{M}(\mathbf{r}_n+i_xL_x\hat{x}+i_yL_y\hat{y}+i_zL_z\hat{z})=\mathbf{M}(\mathbf{r}_n).
\end{align}
\end{subequations}
The integer indices $i_x$, $i_y$ and $i_z\in [-\infty, \infty]$ if the periodicity is present at that specific direction, representing the image's order with respect to the original $0^{th}$ unit cell ($i_x=i_y=i_z=0$).

Cases of T-PBC (Fig.~\ref{fig:t-np} and Fig.~\ref{fig:t-p}) need to be given a special care of as we mentioned before. In these T-PBC cases, the $\mathbf{M}_j$, as a periodic image of $\mathbf{M}_i$ on the PBC, is defined by Eq.~\eqref{eq:pm} that solely depends on the value of its periodic counterpart $\mathbf{M}_i$, thus reducing the numerical problem size to less than $N$. The $(i,j)$ index pairs that satisfy these requirements are referred to as PBC pairs, and we also refer the $i^{th}$ point as the parent of the $j^{th}$ point, and $j^{th}$ point as the child of the $i^{th}$ point. Furthermore, we merge these pairs to their lowest common ancestors (LCA) by the union-find algorithm \cite{10.1145/321879.321884}. Namely, for two PBC pairs $(i,j)$ and $(j,k)$, they are merged to $(i,j)$ and $(i,k)$. These can be proven by adding up their coordinates. A single LCA can have multipe children, e.g., a corner point of a cubic unit cell with 3D PBC can have up to $7$ children. Hence, the number of all unique LCAs $N_{LCA}$ is equal or less than the number of PBC pairs $N_{PBC}$, i.e. $N_{LCA} \leq N_{PBC}$. For a problem of size $N$ with $N_{PBC}$ PBC pairs, only $N-N_{PBC}$ points are unknown, reducing the actual problem size down to $N'=N-N_{PBC}$, so that $N' \leq N$. For simplicity, we assume that the first $N'$ points are within the $0^{th}$ unit cell, and the rest $N-N'$ points are images due to the PBC. 

Note that T-PBC can only appear at the direction with periodicity, yet this is not always true when periodicity is present. For example, consider the 3D periodic case with only $x-$ and $y-$ direction having T-PBC while the $z-$direction having a NT-PBC, e.g., an infinite set of infinite $x-y$ plane films of thickness $d$ stacking over the $z-$direction, with a gap of $L_z-d$.


The discretized LLG equation is written similar to its continuous counterpart:
\begin{equation} \label{eq:llg}
\frac{\partial \mathbf{M}_n}{\partial t} = -\frac{\gamma_{\text{gr}}}{1+\alpha^2}\left(\mathbf{M}_n\times \mathbf{H}_\text{eff,n} + \frac{\alpha}{M_{s\text{,n}}}\mathbf{M}\times\mathbf{M}_n\times\mathbf{H}_\text{eff,n}\right),
\end{equation}
where the solutions are obtained for $\mathbf{M}_n$ and $\mathbf{H}_{\text{eff},n}$ are obtained at the locations $\mathbf{r}_n$ based on the FEM representation of $\mathbf{M}$ via the basis functions in Eq.~\eqref{eq:M_phi}. The time integration can be accomplished via various techniques\cite{Argyris_Scharpf_1969, https://doi.org/10.1002/nme.729}. We implemented implicit time step adaptive predictor-corrector schemes based on explicit second order Adams-Moulton (midpoint) approach as well as time step and order adaptive backward differentiation formulas\cite{1993, Hairer1996}, both approaches having absolute stability up to the second-order time integration.

\section{Implementation} \label{sec:imple}

In this section, we introduce methods to handle periodicity for computing the effective magnetic field components. The anisotropy field $\mathbf{H}_\text{an}(\mathbf{r}_n)$ and applied field $\mathbf{H}_\text{ap}(\mathbf{r}_n)$ are local in that they are solely determined by the their location, so that their computation is straightforward. Approaches for computing the exchange and magnetistatic fields are described next. 

For exchange interaction, modification is only required when PBC present, i.e. $N' < N$. For long-range magnetostatic interaction, the periodic scalar potential (PSP) is via a modified pFFT which employs periodic Green's function as the kernel\cite{AI2024109291}. This can later be further accelerated by reducing the FFT size due to the periodicity of the kernel [ref].

\subsection{Exchange interaction} \label{ssec:ex}

The exchange interactions are due to quantum mechanic effects that originates from Pauli repulsion [ref]. In the continuous micromagnetic representation, the exchange field is given by
\begin{equation} \label{eq:hex}
    \mathbf{H}_{ex}(\mathbf{r}) = \frac{2A_{ex}}{M_s^2(\mathbf{r})}\nabla^2 \mathbf{M}(\mathbf{r}).
\end{equation}
where $\nabla^2$ is the Laplacian operator. The numerical representation of the exchange field is given by
\begin{equation} \label{eq:lapl}
    \mathbf{H}_{ex}(\mathbf{r}_n)=\sum_{\mathbf{r}_m\in \Omega_n}\omega_{mn}\mathbf{M}(\mathbf{r}_{m})
\end{equation}
where $n\in[1,N']$, $m\in[1,N]$, $\omega_{mn}$ are elements of a sparse matrix that has band of non-zero values defined via $\Omega_n$, such that each point is connected to the surrounding points through their connecting mesh edges\cite{Fastmag}. For the NP-PBC case, the calculation of the exchange matrix is identical to that of the conventional case of non-periodic problems. 

For the case of T-PBC, the sparse matrix construction needs to be modified for several cases.

\textit{Case 1:} This case is for the interior points, such that the point at $\mathbf{r}_n$ is not the LCA of any other points and all the neighboring points within $\Omega_n$ are not LCA of any other points, i.e., $LCA(i)\neq n,i\in[1,N],i\neq n$ and $\forall m\in\Omega_n$ we have $LCA(j)\neq m,j\in[1,N],j\neq m$. Namely, both $n$ and $m\in\Omega_n \in [1,N']$. This case represents the majority of the sparse matrix construction since most points within the unit cell volume that are not close to a PBC boundary fall within this category. The resulting sparse matrix elements of Eq.~\eqref{eq:lapl} are the same as for the non-periodic case.

\textit{Case 2:} In this case, the point at $\mathbf{r}_n$ is not the LCA of any other points, and some of the points within $\Omega_n$ have different LCAs other than themselves. We need to map these points at $\mathbf{r}_m$ to their LCAs and Eq.~\eqref{eq:lapl} can be extended as
\begin{equation} \label{eq:lapl_1}
    \mathbf{H}_{ex}(\mathbf{r}_n)=\sum_{\mathbf{r}_m\in \Omega_n}\omega_{mn}\mathbf{M}(\mathbf{r}_{LCA(m)})
\end{equation}
to retrieve the magnetization on the PBC from within $0^{th}$ unit cell. The weight $\omega_{mn}$ is not mapped, since though the magnetization is from the $\mathbf{r}_{LCA(m)}$, the location should still stay within $\Omega_n$. Note that if $m$ is not a child in a PBC pair, then its $LCA$ is itself, i.e., $LCA(m)=m,m\in[1,N']$, which falls back to the original form of Eq.~\eqref{eq:lapl}.

\textit{Case 3:} In this case, the point at $\mathbf{r}_n$ is the LCA of some points on the boundary, i.e., $\exists n', n'\in[N'+1,N]$ and $LCA(n')=n$. Following the same approach, for each $n'$ we can define its surrounding region as $\Omega_{n'}$. Then, the completed surrounding region of $\mathbf{r}_n$ is $\Omega_{n}\cup_{n'}\Omega_{n'}$, denoted as $\Omega_{n\cup \{n'\}}$. Here, $\{n'\}$ represent the union set of all $n'$. We can see that such new region falls back to the original $\Omega_n$ if $\{n'\}=\varnothing$. Similar to Eq.~\eqref{eq:lapl}, we define all the points inside the $\Omega_{n\cup \{n'\}}$ with index $m$. We further note that an edge can occur multiple times in this object of size $N$, hence direct computations via Eq.~\eqref{eq:lapl_1} would introduce duplicated calculations. To avoid such duplication, we first defined a new set of index $m'$ and corresponding coordinates $\mathbf{r}_{m'}$. For any $m\in\Omega_{n\cup \{n'\}}$, a new $m'$ and $\mathbf{r}_{m'}$ can be generated as
\begin{equation} \label{eq:shift}
    \begin{aligned}
        &LCA(m')=LCA(m) \\
        &\mathbf{r}_{m'} - \mathbf{r}_m = \mathbf{r}_{n'} - \mathbf{r}_n
    \end{aligned}
\end{equation}
Combining all $\mathbf{r}_{m'}$ and $m'$ together, we obtain a new location set $\{\mathbf{r}_{m'}\}$ and the corresponding indices set $\{m'\}$. Then, the desired locations and indices set represent a unique set of all $\mathbf{r}_m'$, written as $\{\mathbf{r}_m'\}^\dagger$ and their corresponding indices, denoted as $\{m'\}^{\dagger}$. Here, $\{\mathbf{r}_m'\}^\dagger$ and $\{m'\}^{\dagger}$ satisfies
\begin{equation} \label{eq:dedup}
    \begin{aligned}
        &\{\mathbf{r}_{m'}\}^\dagger\subseteq \{\mathbf{r}_m'\},\{m'\}^{\dagger} \subseteq \{m'\}\\\
        &\forall i\in \{m'\},\exists j\in\{m'\}^{\dagger}\Rightarrow \mathbf{r}_i=\mathbf{r}_j ~~\mathbf{r}_i,\mathbf{r}_j\in\{\mathbf{r}_{m'}\}^\dagger\\
        &\forall \alpha,\beta\in\{m'\}^\dagger,\alpha\neq\beta\Rightarrow\mathbf{r}_\alpha\neq\mathbf{r}_\beta
    \end{aligned}
\end{equation}
This new location set can, then, form a new region as $\Omega_{m'}^\dagger$, which, combined with the original region $\Omega_n$, determines a generalized surrounding region defined as $\Tilde{\Omega}_n=\Omega_n\cup\Omega_{m'}^\dagger$. Finally with this new region, we can rewrite the Eq.~\eqref{eq:lapl_1} as
\begin{equation} \label{eq:lapl_2}
    \mathbf{H}_{ex}(\mathbf{r}_n)=\sum_{\mathbf{r}_{k}\in\Tilde{\Omega}_n}\omega_{kn}\mathbf{M}(\mathbf{r}_{LCA(k)})
\end{equation}

A fast implementation to construct the above generalized surrounding region defined in Eq.~\eqref{eq:dedup} is to place a filter in Eq.~\eqref{eq:shift}. Such filter rejects the indices between $N'+1$ and $N$, i.e., we only map the non-PBC location. In this case, the equality $\{\mathbf{r}_{m'}\}^\dagger= \{\mathbf{r}_m'\},\{m'\}^{\dagger} = \{m'\}$ holds directly.

By combining Eq.~\eqref{eq:hex} and Eq.~\eqref{eq:lapl}, Eq.~\eqref{eq:lapl_1}, Eq.~\eqref{eq:lapl_2}, we can handle the exchange field when T-PBC is present.

\subsection{Magnetostatic interaction} \label{ssec:ms}
\subsubsection{Problem formulation}
The magnetostatic field is due to long-range interactions, and it can be calculated either by solving the Poisson equation or by evaluating the superposition integrals \cite{NUFFT, 10259661}. Here, we evaluate the magnetostatic field using superposition integrals via the following formulation
\begin{subequations} \label{eq:hms}
    \begin{align}
        & \rho(\mathbf{r}) = \nabla\cdot\mathbf{M}(\mathbf{r}),~\rho_s(\mathbf{r})=-\hat{\mathbf{n}}\cdot\mathbf{M}(\mathbf{r}),\label{eq:hms1}\\
        & u(\mathbf{r})=\iiint_{V}G^p(\mathbf{\mathbf{r}-\mathbf{r'}})\rho(\mathbf{r})d\mathbf{r'}+\iint_{S}G^p(\mathbf{\mathbf{r}-\mathbf{r'}})\rho_s   (\mathbf{r})d\mathbf{r'},\label{eq:hms2}\\
        & \mathbf{H}_{ms}(\mathbf{r}) = -\nabla u(\mathbf{r}).\label{eq:hms3}
    \end{align}
\end{subequations}
Here, $\rho$ is the volumetric magnetic charge density, $\rho_s$ is the surface magnetic charge density, $u$ is the magnetic scalar potential, and $G_p$ is the 3D periodic Green's function (PGF), which can be 1D, 2D, or 3D periodic \cite{AI2024109291}. A discrete FEM representation of Eq.~\eqref{eq:hms} can be given as \cite{Fastmag, NUFFT}
\begin{subequations} \label{eq:hms_d}
    \begin{align}          &q(\mathbf{r}_n)=\sum_{\mathbf{r}_k\in\Tilde{\Omega}_n}\omega^{q}_{kn}\mathbf{M}(\mathbf{r}_{LCA(k)}),\label{eq:hms_d1}\\
        & \tilde{u}(\mathbf{r}_n) = \sum_{k\in[1,N']} G^p(\mathbf{r}_n-\mathbf{r}_k)q(\mathbf{r}_k),\label{eq:hms_d2}\\
        & u(\mathbf{r}_n) = \tilde{u}(\mathbf{r}_n) + \sum_{\mathbf{r}_k\in\Omega^{\text{near}}_n} \omega^{\text{near}}_{kn} \mathbf{M}(\mathbf{r}_k),\label{eq:hms_d2b}\\
        & \mathbf{H}_{ms}(\mathbf{r}_n)=-\sum_{\mathbf{r}_k\in\Tilde{\Omega}_n}\omega^\text{grad}_{kn}u(\mathbf{r}_{LCA(k)}).\label{eq:hms_d3}
    \end{align}
\end{subequations}
Here, $q(\mathbf{r}_n)$ are nodal magnetic charges assigned to each point $\mathbf{r}_n$, i.e., nodes of the FEM mesh. These charges are found via an FEM representation of the divergence operator in the volume and additional surface related charge components on the surfaces of the magnetization discontinuities, e.g., boundaries of the magnetized domain or boundaries between material of different saturation magnetization. The nodal charges are given as a sparse matrix-vector products with the sparse matrix elements $\omega^{q}_{kn}$. Similarly, $\omega^\text{grad}_{kn}$ are elements of the sparse matrix representing the gradient operator. For the NT-PBC cases, the charge and gradient sparse matrices are constructed identically to non-periodic cases. For the T-PBC cases, the construction of these sparse matrices is modified similarly to the modifications in the Laplace matrix construction, as outlined in Sec.~\ref{ssec:ex} to generalize the regions surrounding each point. 

In Eq.~\eqref{eq:hms_d2}, $\tilde{u}(\mathbf{r}_n)$ is the magnetic periodic scalar potential (PSP) at the mesh points and it is given as a standard superposition sum over the charges $q(\mathbf{r}_k)$ in the periodic unit cell, where $G^p(\mathbf{r}_n-\mathbf{r}_k)$ is the periodic Green's function (PGF) defined for the 1D, 2D, and 3D periodicities as infinite sums \cite{FPIP}:
\begin{subequations}
    \label{eq:2}
    \begin{align}
        G^{p}(\mathbf{r}) =& \sum_{i_x=-\infty}^{\infty} e^{-jk_{x0}i_xL_x}G_0(\mathbf{r}-i_xL_x\mathbf{\hat{x}}) \text{     for 1D},\\
        G^{p}(\mathbf{r}) =& \sum_{i_x=-\infty}^{\infty}\sum_{i_y=-\infty}^{\infty}e^{-j[k_{x0}i_xL_x+k_{y0}i_yL_y]}G_0(\mathbf{r}-i_xL_x\mathbf{\hat{x}}-i_yL_y\mathbf{\hat{y}}) \text{     for 2D},\\
        G^{p}(\mathbf{r}) =& \sum_{i_x=-\infty}^{\infty}\sum_{i_y=-\infty}^{\infty}\sum_{i_z=-\infty}^{\infty}e^{-j[k_{x0}i_xL_x+k_{y0}i_yL_y+k_{z0}i_zL_z]}G_0(\mathbf{r}-i_xL_x\mathbf{\hat{x}}-i_yL_y\mathbf{\hat{y}}-i_zL_z\mathbf{\hat{z}})  \text{ for 3D},
    \end{align}
\end{subequations}
where $G_0(\mathbf{r})=1/|\mathbf{r}|$ is the free-space Green's function. The PSP can be calculated efficiently \cite{AI2024109291} for the NP-PBC case and these approaches can be adapted to account for the P-PBC cases as discussed next. 

For a given charge used for the calculation of the PSP, the PSP values are calculated accurately only at distances substantially distant from the charge, specifically, for distances greater than the surrounding element size. For smaller distances, the result is inaccurate and it needs to be corrected. In Eq.~\eqref{eq:hms_d2b}, $u(\mathbf{r}_n)$ is the corrected PSP (C-PSP), which includes a correction term, where $\omega^{\text{near}}_{kn}$ are elements of the sparse correction matrix. The correction matrix has non-zero elements for a range of points $\Omega^{\text{near}}_n$ in the proper vicinity of $\mathbf{r}_n$. The correction sparse matrix elements are found through exact integrals of Eq.~\eqref{eq:hms_d2} performed over the relevant elements (tetrahedrons in our case). Approaches to account for the PBCs in the computations of Eq.~\eqref{eq:hms_d} are outline next. With such a correction, the C-PSP can be computed with a controllable error.

\subsubsection{Implementation}
 Equation~\eqref{eq:hms_d2} includes a superposition sum and if computed directly it has its computation cost scaling as $O(N^2)$, which is high for large $N$. We use a modified version of the pFFT method, also referred to as box-adaptive integral method (BAIM)\cite{https://doi.org/10.1029/96RS02504} that has a computational complexity of $O(N\log N)$. We define a Cartesian grid of size $N_x\times N_y\times N_z$ over the domain size of $D_x, D_y, D_z$ as
\begin{equation}\label{eq:grid}
\begin{aligned}
    &x_{i,j,k} = (i-1) \Delta_x, \quad y_{i,j,k} = (j-1)\Delta_x, \quad z_{i,j,k} = (k-1)\Delta_x\\
    &\Delta_x = \frac{D_x}{N_x-1}\quad \Delta_y = \frac{D_y}{N_y-1}\quad \Delta_z = \frac{D_z}{N_z-1}
\end{aligned}
\end{equation}
with $i\in[1,N_x],j\in[1,N_y],k\in[1,N_z]$ and denote the grid points as $\mathbf{r}_{i,j,k}=(x_{i,j,k},y_{i,j,k},z_{i,j,k})$. The BAIM procedure follows the following four steps.

\textit{step 1:} Project non-uniformly distributed sources $q(\mathbf{r}_n)$ of Eq.~\eqref{eq:hms1} from locations $\mathbf{r}_n$ to the charges $q_{i,j,k}$ defined at the uniform Cartesian grid at $\mathbf{r}_{i,j,k}$.

\textit{step 2:} Calculate the PSP on all Cartesian grid points $\mathbf{r}_{i,j,k}$ as a convolution between PGF on  $G_p(\mathbf{r}_{i,j,k})$ and changes $q_{i,j,k}$ on the grid, obtaining the PSP $u_{i,j,k}$ on the grid.

\textit{step 3:} Interpolate the grid PSP $u_{i,j,k}$ to the PSP at the non-uniform points $\mathbf{r}_n$.

The result from Steps 1-3 provide a PSP at the non-uniform points that has errors due to the fact that the projection and interpolation procedures are inaccurate from the PSP generated by sources that are located near the observers. To correct these inaccuracies, we implement a correction Step 4. 

\textit{step 4:} Correct the errors introduced by projection and interpolation from Steps 1 and 3 to obtain the desired PSP $u(\mathbf{r}_n)$ with a controllable error level.

Compared with the non-periodic problems, two major modifications are made in these four steps. 

First, in the original BAIM, Green's function in step 2 is the free space Green's function $G_0=1/|\mathbf{r}|$. To adapt the PBCs, the free-space Green's function is replaced with the PGF $G_p$ and it can be computed efficiently \cite{LGF, AI2024109291}.

The second major modification is for the error correction in Step 4. In general, the projection and interpolation only introduce large errors when projected/interpolated functions vary drastically within the region of interest. For the free-space Green's function $G_0$ it occurs when the separation between source and observer is small. We define a threshold $|\mathbf{r}_{ER}|$ such that for all neighboring sources within this range, contributions from these sources are calculated directly instead of the convolution in step 2. This can be done by subtracting the contributions from step 2 and adding the accurate contributions computed  via a direct contribution. 

For the free-space case, the threshold $|\mathbf{r}_{ER}|$ is found simply based on the source-observer distance. For the PBC cases, the periodic source images result in additional periodically copied points that also need to be included in the threshold $|\mathbf{r}_{ER}|$ range. Although we replace $G_0$ with $G_p$ in Step 2, only $G_0$ varies drastically within nearby region [ref]. Therefore regardless of the periodicity configuration, the error correction is always calculating the contributions from sources through $G_0$, namely we only subtract contributions of $G_0$ in step 2 and add them back still using $G_0$.

The $|\mathbf{r}_{ER}|$ here defines a correction region for each point as $\Omega_n^{ER}$ that includes all qualified neighboring source points, this is usually defined based on the box size $\Delta_x,\Delta_y,\Delta_z$, e.g. $|\mathbf{r}_{ER}|=max(\Delta_x,\Delta_y,\Delta_z)$. This region need to be expanded when periodicity presents. Similar to Sec.~\ref{ssec:ex} we have different scenarios. We first define a symbol $\iota$ to represent the periodic direction, e.g. $\iota$ can be $x$ for 1D periodicity, $x,y$ for 2D periodicity and $x,y,z$ for 3D periodicity. Then for a specific point at location $\mathbf{r}_n$ we have

\textit{case 1:} $\forall \iota\Rightarrow\mathbf{r}_n\hat{\iota} \pm |\mathbf{r}_{ER}| \in (D_{\iota}-L_{\iota},L_{\iota})$. This is the major case for all source points that are sufficiently far away from the boundary of periodic direction. In this case $\Omega_n^{ER}$ is unchanged.

\textit{case 2:} $\exists \iota\Rightarrow \mathbf{r}_n\hat{\iota}+|\mathbf{r}_{ER}| \geq L_{\iota}$ or $\mathbf{r}_n\hat{\iota}-|\mathbf{r}_{ER}| \leq D_{\iota}-L_{\iota}$. This means that the error correction region is including source points from nearest neighboring images. The error correction region then need to be expanded. The original $\Omega_n^{ER}$ is defined within a region that satisfies $\forall \mathbf{r}_k \in \Omega_n^{ER} \Rightarrow \forall \iota,\mathbf{r}_k\hat{\iota} \in [\mathbf{r}_n\hat{\iota}-|\mathbf{r}_{ER}|,\mathbf{r}_n\hat{\iota}+|\mathbf{r}_{ER}|]$. We here defined a expanded error correction as $\Tilde{\Omega}_n^{ER}$. For each scenario we have
\begin{equation} \label{eq:new_er}
    \begin{aligned}
        &\forall \mathbf{r}_k \in \Tilde{\Omega}_n^{ER}\Rightarrow \forall \iota,\mathbf{r}_k\hat{\iota} \in [\mathbf{r}_n\hat{\iota}-|\mathbf{r}_{ER}|,\mathbf{r}_n\hat{\iota}+|\mathbf{r}_{ER}|] \cup [0, \mathbf{r}_n\hat{\iota}+|\mathbf{r}_{ER}| - L_{\iota}]\\
        &\forall \mathbf{r}_k \in \Tilde{\Omega}_n^{ER}\Rightarrow \forall \iota,\mathbf{r}_k\hat{\iota} \in [\mathbf{r}_n\hat{\iota}-|\mathbf{r}_{ER}|,\mathbf{r}_n\hat{\iota}+|\mathbf{r}_{ER}|] \cup [L_{\iota}+\mathbf{r}_n\hat{\iota}-|\mathbf{r}_{ER}|,D_{\iota}]
    \end{aligned}
\end{equation}
In BAIM, the source points within each range can be easily found by searching "boxes" constructed by grid points of size $\Delta_x,\Delta_y,\Delta_z$. And all these source points' corresponding box index are pre-tabulated.

The construction of the sparse matrix $\mathbf{Z_0}(\mathbf{r}_n)$ in Eq.~\eqref{eq:hms_d3} also involves very similar modification, except that the error correction range $|\mathbf{r}_{ER}|$ is typically defined based on the tetrahedron size instead of the box size, i.e. the average edge length. Otherwise we also need to expand the region searching for neighbor following the same style defined in in Eq.~\eqref{eq:new_er}.

Here we first assume that the case is NP-PBC. This assumption naturally holds for most regular shapes, e.g. a cube, a cylinder and etc. However, this property is not guaranteed to be held with an irregular shape. In this case the judging condition in case 1 and 2 are not working anymore.

To fix this issues, here we propose a very simple solution. When calculating the $H_{ms}$ over the unit cell that $\exists \iota\Rightarrow D_{\iota} > L_{\iota}$. We define a new unit cell based on the original first $N'$ points specifically for $H_{ms}$ that
\begin{equation}\label{eq:shift_hms}
    \begin{aligned}
        &\mathbf{r}_{n'} =  \underset{\delta} {\mathrm{argmin}}\sum_{\iota} |\mathbf{r}_n\hat{\iota}-\mathbf{r}_0\hat{\iota}+\delta L_{\iota}|,\delta\in \mathbb{Z}\\
        &H_{ms}(\mathbf{r}_{n'}) = H_{ms}(\mathbf{r}_n)
    \end{aligned}
\end{equation}
This equation defines a new unit cell with new coordinates $\mathbf{r}_{n'}$. The coordinates are shifted by arbitrary periodic length along each periodic direction to minimize the Manhattan distance with respect to a fixed point $\mathbf{r}_0$. For convenience we typically choose the geometry center of the original unit cell, namely $\mathbf{r}_0=\sum_n\mathbf{r}_n/N'$. The new unit cell defined by Eq.~\eqref{eq:shift_hms} satisfies the condition that $\forall \iota \Rightarrow D_{\iota} \leq L_{\iota}$. Since $\forall \mathbf{r}_{n'}\hat{\iota}\Rightarrow \mod(\mathbf{r}_n\hat{\iota}-\mathbf{r}_0\hat{\iota},L_{\iota})\in[0,L_{\iota})$. After solving the $H_{ms}$ on each shifted point, we can map all the values from $\mathbf{r}_{n'}$ back to $\mathbf{r}_n$ to restore the original unit cell.

\section{Results} \label{sec:res}

PM-FEM was implemented as a part of the FastMag micromagnetic simulator to allow modeling of complex magnetic materials and devices\cite{Fastmag}. The updated FastMag components include the magnetostatic and exchange field modules as well as the mesh construction module for touching and protruding unit cell types. In addition to these modules, FastMag implements implicite time stepping schemes, including time step and order adaptive backward diffrentiation formulas and time step adaptive second-order midpoint rule. In these implicit time stepping schemes, at each time step, a non-linear equation is solved via iterative Newton's method assisted with a linear solver at each non-linear iteration. The convergence of the linear solver is assisted with a preconditioner, such as ILU or block inverse preconditioners\cite{Scott2023-th, 8902292}. 

We first show a transformation of a 1D periodic protruding touching unit cell as in \ref{ssec:prot}. We then present examples of micromagnetic simulations for 1D/2D/3D periodic problems. Finally, we show the computational performance of PM-FEM. We do not show results for computing the periodized exchange field since the computational time for the periodic and non-periodic cases of computing the exchange field is nearly the same and this time is significantly below that of themagnetostatic field.

\subsection{Protruding unit cell} \label{ssec:prot}

\begin{figure}[htbp] \label{fig:prot}
    \centering
    \subfigure[]{
    \includegraphics[width=0.28\linewidth]{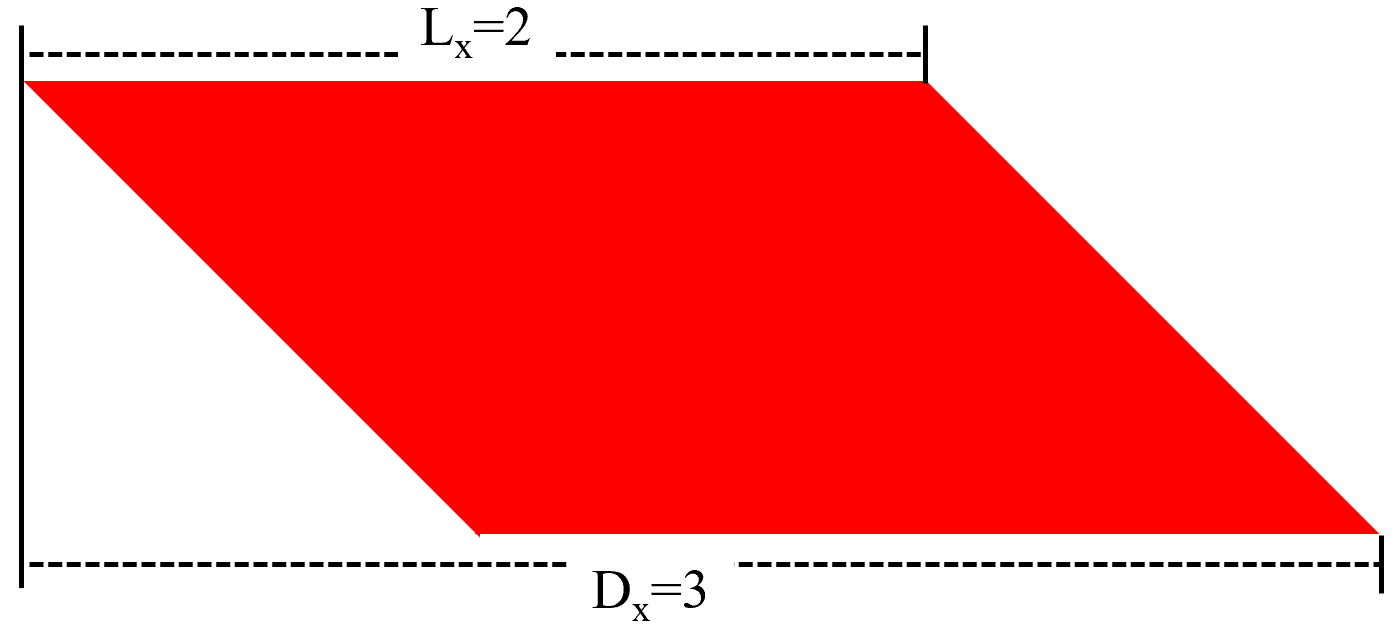} \label{fig:prot_1}
    }
    \quad
    \subfigure[]{
    \includegraphics[width=0.29\textwidth]{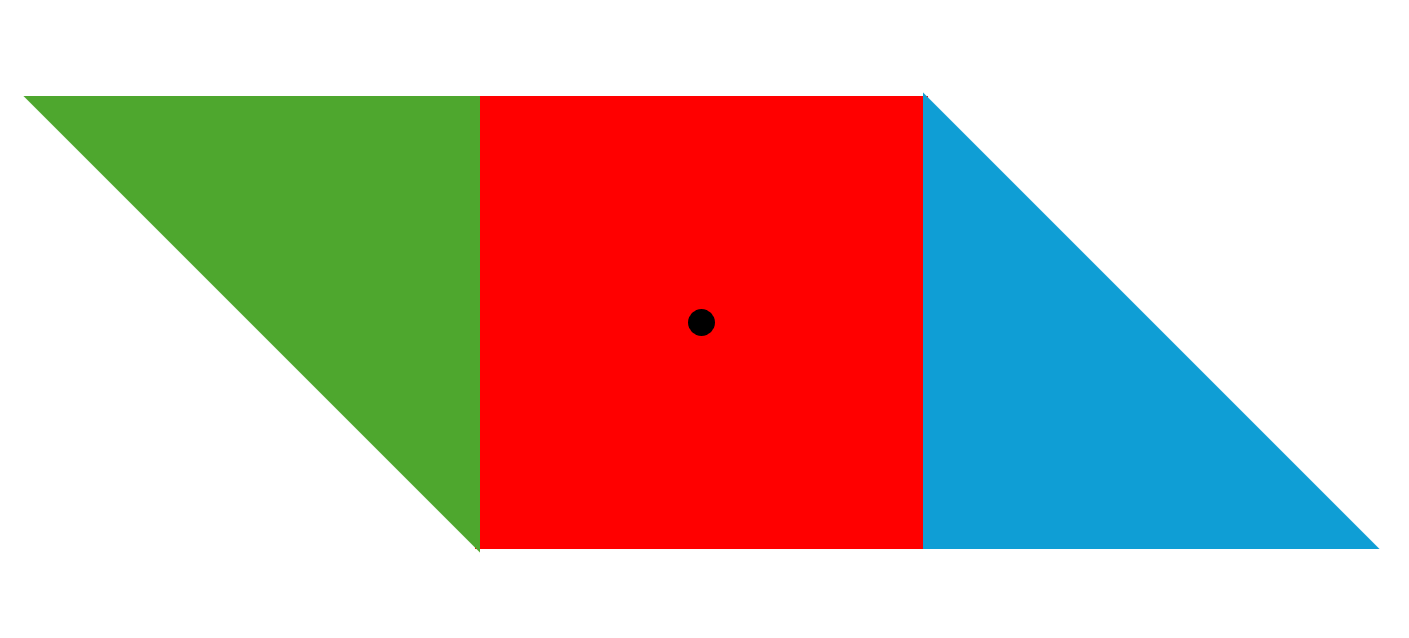} \label{fig:prot_2}
    }
    \quad
    \subfigure[]{
    \includegraphics[width=0.2\textwidth]{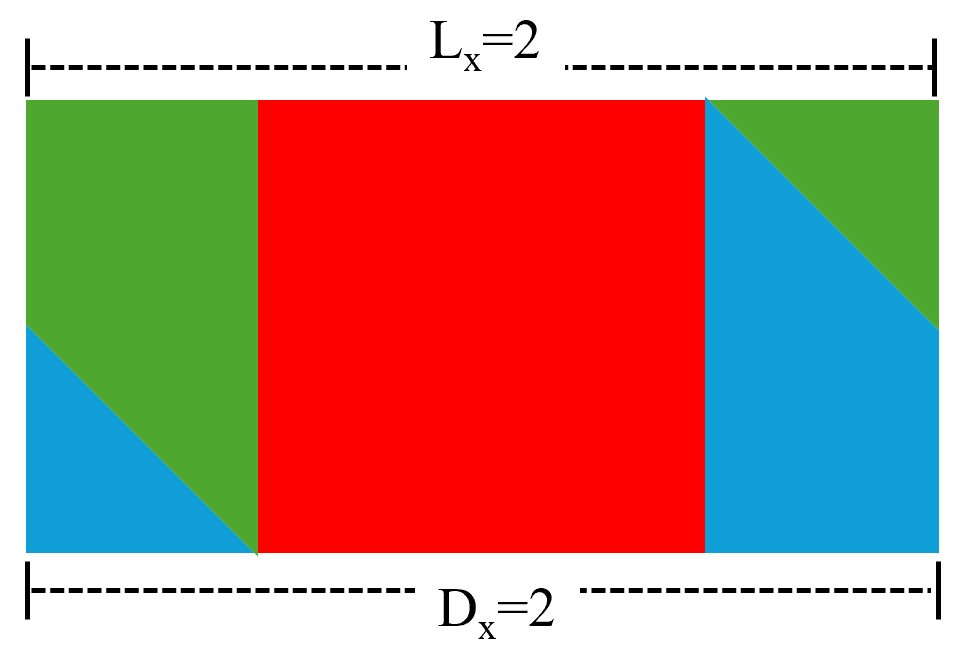} \label{fig:prot_3}
    }
    \caption{Illustration of (a) Protruding unit cell; (b) its geometry center (black dot) and its protruding parts in green and blue; (c) Regular unit cell after shifting the protruding parts.}
    \label{fig:shift}
\end{figure}

The protruding unit cell is defined such that there is at least one periodic direction where the structure size exceeds the periodic length. In Fig.~\ref{fig:prot_1}), the structure is a parallelogram with a 1D $x-$direction T-P-PBC. The cell size along the $x-$direction is one and a half of the periodic length, e.g., $L_x=2,D_x=3$. By utilizing Eq.~\eqref{eq:shift_hms} with $\mathbf{r}_0$ set to the geometry center of the unit cell (the black dot in Fig.~\ref{fig:prot_2}), we shift the green and blue part in Fig.\ref{fig:prot_2} via Eq.~\eqref{eq:shift_hms}, and obtain a new unit cell with a regular shape, i.e., a rectangular unit cell with T-NP-PBC (Fig.~\ref{fig:prot_3}). This step is only needed once in the preprocessing stage to map the coordinates, and it has no effect on the rest of the computations.

\subsection{1D periodicity} \label{ssec:1d}

\begin{figure}
    \centering
    \includegraphics[width=0.6\linewidth]{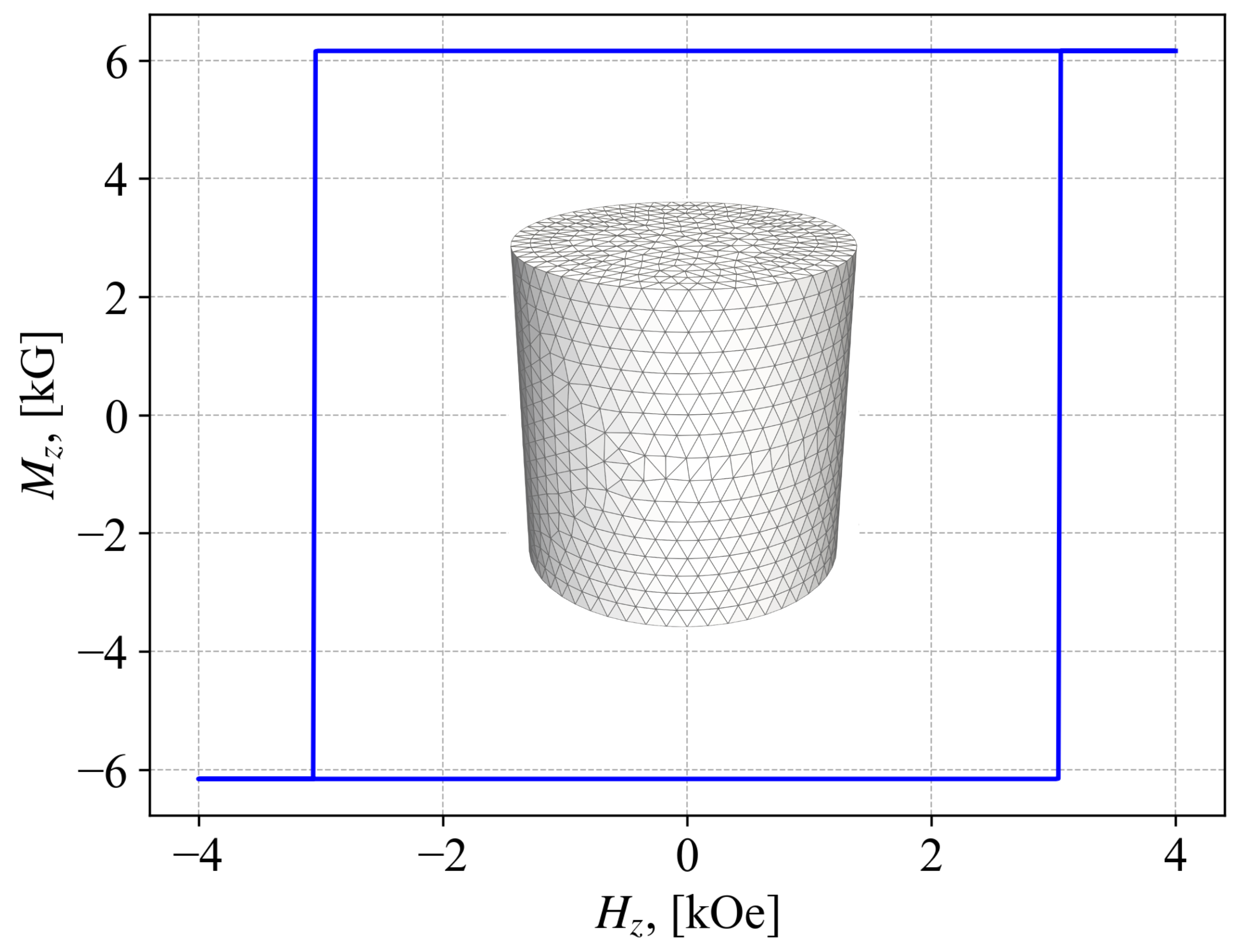}
    \caption{Hysteresis loop along $\hat{z}-$direction of infinite long periodic rod along $\hat{z}-$direction, the inset is the unit cell. The loop is of square shape and coercive field is around 3050 Oe, close to theoretical value 3085 Oe.}
    \label{fig:1d}
\end{figure}

Figure~\ref{fig:1d} shows a calculated M-H hysteresis loop for an infinitely long cylindrical rod with $A_{ex}=9.604\times 10^{-7}~\text{erg/cm}$, $M_s=490~\mathrm{emu/cm^3}$, $\alpha=0.5$, and no anisotropy. The corresponding exchange length is $\lambda_{ex}=\sqrt{A_{ex}/2\pi M_s^2}\approx 8\times10^{-7}~\text{cm}$. The prediction \cite{PhysRev.109.1522} is a mostly square hysteresis loop and T-PBC is required especially for the case where the radius of the rod $R$ is close to the normalized length $R=R_0=\sqrt{2\pi}\lambda_{ex}\approx2\times10^{-6}~\text{cm}$ due to the buckling reversal mechanism\cite{PhysRev.106.446}. We choose the height of the unit cell twice of the radius, i.e., $h=2R$ (shown as the inset of Fig.~\ref{fig:1d}), as indicated in\cite{1Dtub}, expecting accurate results. We impose a 1D T-NP-PBC along the height ($z-$) direction. The switching is via buckling mode and the coercive field is $H_c== 3085~\text{Oe}$, which agrees with solutions in literature \cite{PhysRev.109.1522}.

\subsection{2D periodicity} \label{ssec:2d}

\begin{figure}[htbp] \label{fig:2D}
    \centering
    \subfigure[]{
    \includegraphics[width=0.35\linewidth]{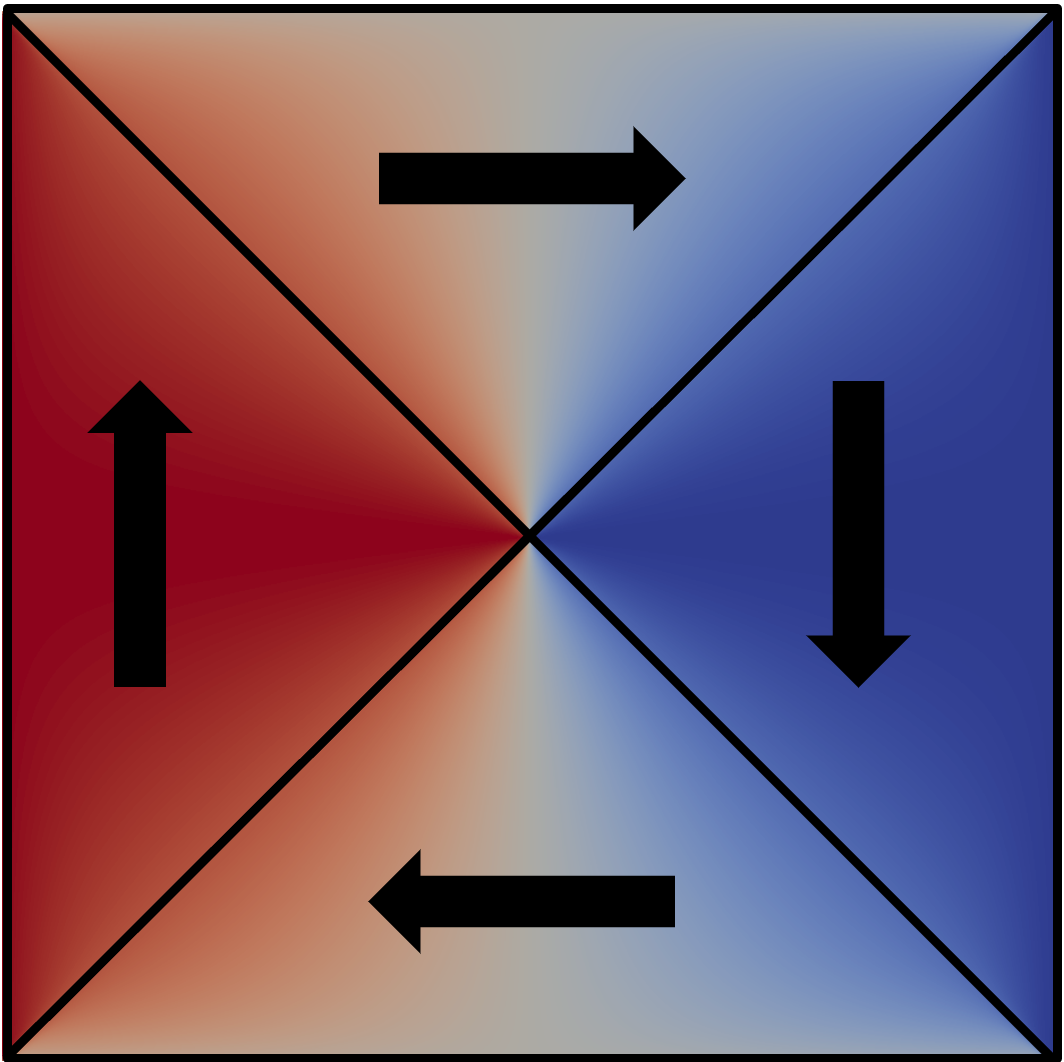} \label{fig:2D_1}
    }
    \quad
    \subfigure[]{
    \includegraphics[width=0.35\textwidth]{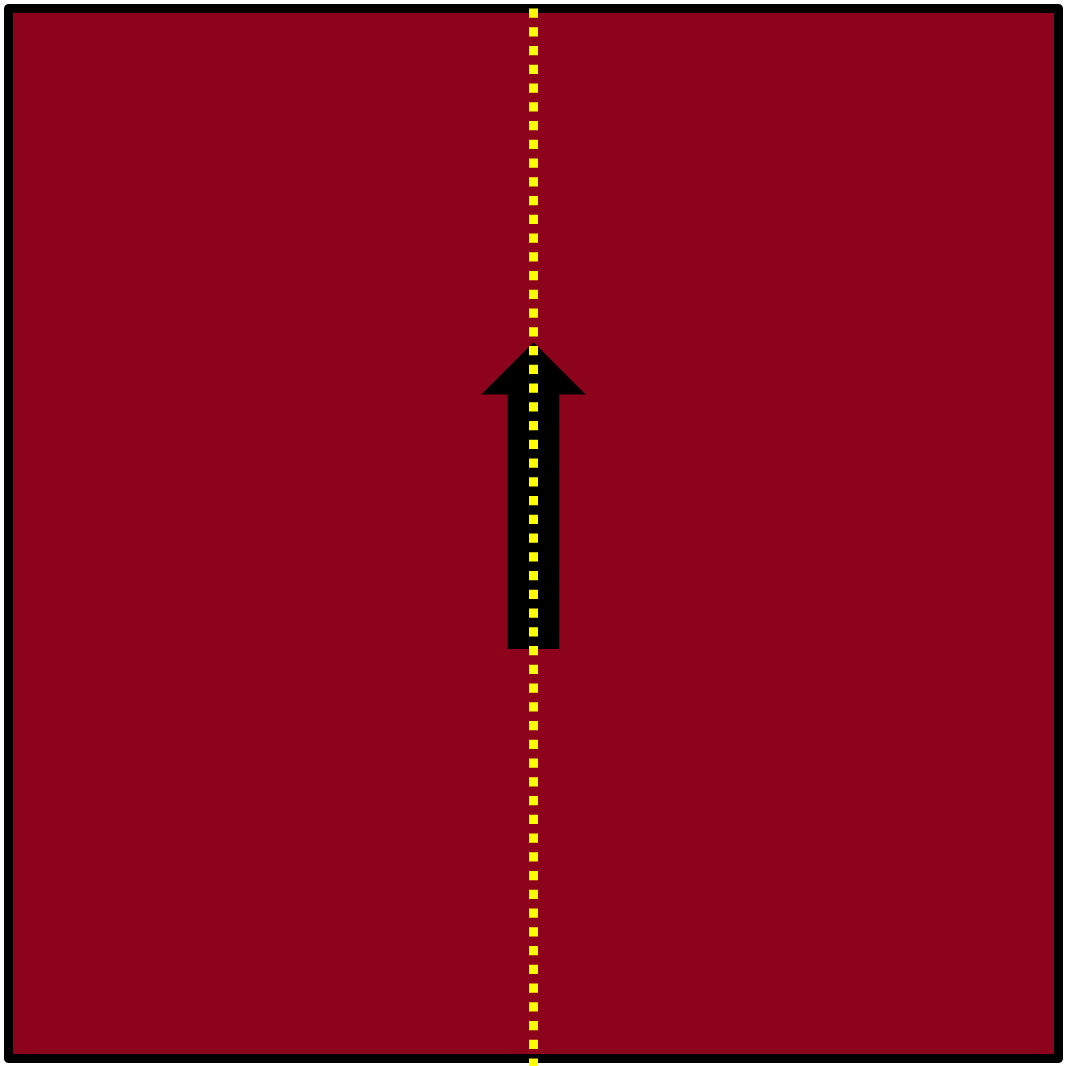} \label{fig:2D_2}
    }
    \newline
    \subfigure[]{
    \includegraphics[width=0.35\textwidth]{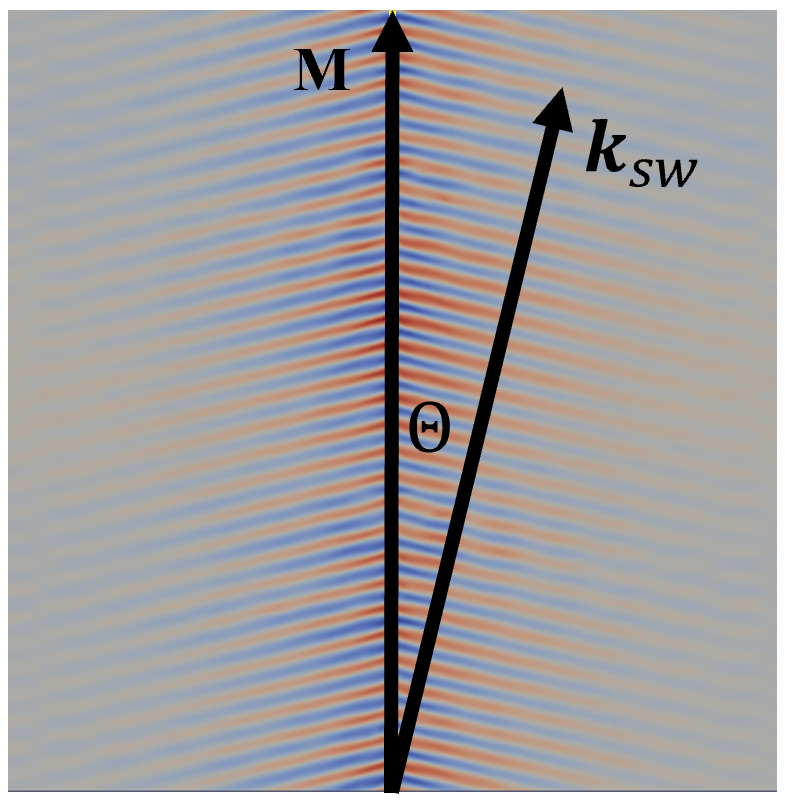} \label{fig:2D_3}
    }
    \quad
    \subfigure[]{
    \includegraphics[width=0.5\textwidth]{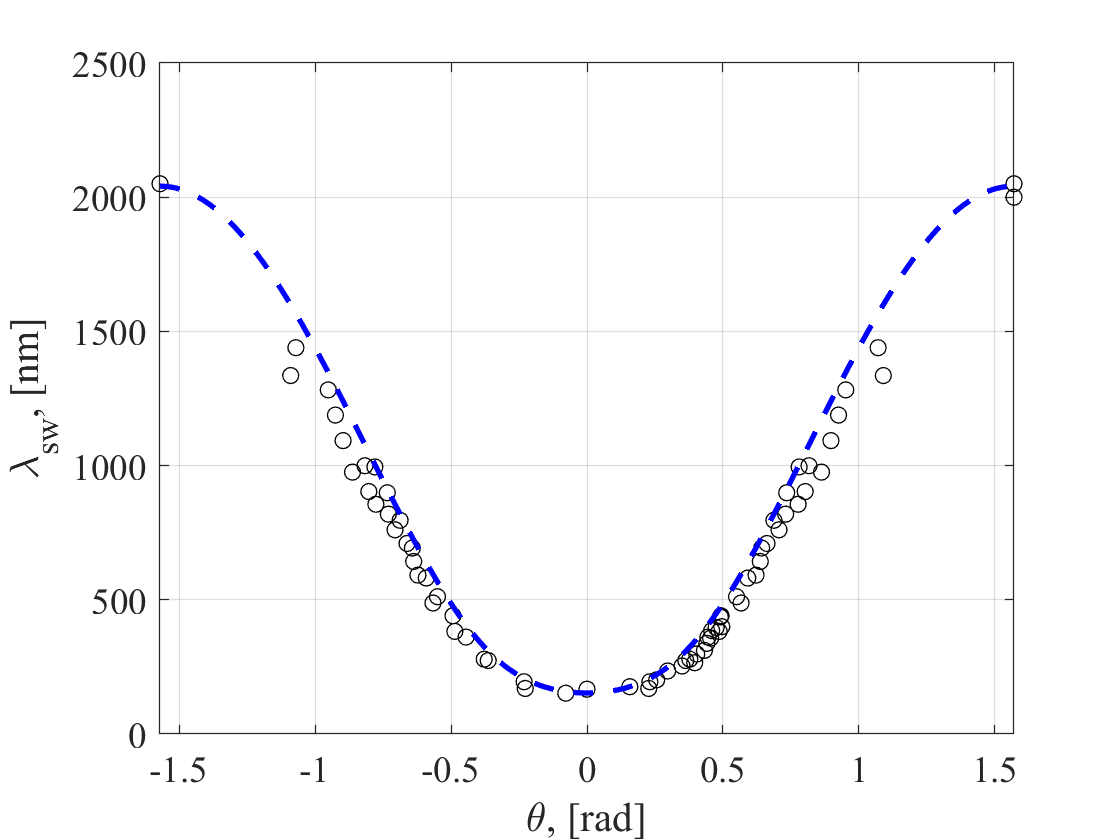} \label{fig:2D_4}
    }
    \caption{Example of 2D T-NP-PBC case. (a) Equilibrium vortex state when no PBC is present and (b) the uniform state when the 2D T-NP-PBC is used. The yellow dashed line is the location of the line excitation. (c) Angles between wave vector of propagating spin wave and magnetization. (d) Simulated wavelength (circle marks) and theoretical prediction (blue dashed line).}
    \label{fig:2d}
\end{figure}

Here, we present results for an infinite Permalloy film, modelled as a 2D T-NP-PBC square unit cell of the side length of  $D_x=D_y=2\times10^{-5}$ cm and thickness $D_z=2\times10^{-6}$ cm (Fig.~\ref{fig:2D_1}). The material parameters are  $M_s=637~\text{emu/cc}$, $A_{ex}=1.4\times 10^{-6}~\text{erg/cm}$, zero anisotropy, and $\alpha=0.02$. We impose a 2D T-PBC in the $\hat{x}$ and $\hat{y}$ directions, i.e., $L_x=L_y=D_x=D_y$. We compare the results for the cases with and without the PBC. Without the PBC, the equilibrium is a vortex state, related to the magnetostatic effects of the edges (Fig.~\ref{fig:2D_1}). With the PBC, the equilibrium state is uniform in the $x-y$ plane (e.g., magnetization aligned along the $y-$direction in Fig.~\ref{fig:2D_2}).

We then use a line source (bright yellow line in Fig.~\ref{fig:2D_2}) in the middle of the film along the $y-$direction parallel to the equilibrium magnetization. The line source represents a radio frequency (RF) excitation $H(\mathbf{r},t)=H_0\cos(\omega t-|\mathbf{k}|y)$ of width $d=10^{-6}$ cm $\ll 2\pi/|\mathbf{k}|$. We fix the driving frequency $\omega_0=4.85\times2\pi$ GRad/s and tune wavelength $2\pi/|\mathbf{k}|$. Due to the phase matching, excited propagating spin waves (Fig.~\ref{fig:2D_3}) share the same wavenumber with tilted angle $\theta$ between the wave vector $\mathbf{k}_{sw}$ and magnetization $\mathbf{M}$ that satisfy $\mathbf{k}_{sw}\cos\theta=\mathbf{k}$. Therefore, by sweeping the $\mathbf{k}$, we can extract the dispersion relationship of the spin wave for the wavelength $\lambda_{sw}=2\pi/|\mathbf{k}_{sw}|$ versus the angle $\theta$. We compare the calculated dispersion relationship with a theoretical result\cite{Kalinikos1986}:
\begin{equation} \label{eq:sw_disp}
    \omega^2(\mathbf{k}_{sw}) = \gamma_{ex}\omega_M\mathbf{k}_{sw}^2 \times \left[  \gamma_{ex}\omega_M\mathbf{k}_{sw}^2 + \omega_M\times\left(1-\frac{|\mathbf{k}_{sw}|D_z\cos^2\theta}{2}  + \omega_M\frac{|\mathbf{k}_{sw}|D_z(2-|\mathbf{k}_{sw}|D_z)\sin^2\theta}{4 \gamma_{ex}\omega_M\mathbf{k}_{sw}^2} \right) \right],
\end{equation}
where $\omega_M=\gamma_\text{gr}M_s$, $\gamma_{ex}=\lambda_{ex}^2=A_{ex}/2\pi M_s^2$. By fixing $\omega(\mathbf{k}_{sw})=\omega_0$ in Eq.~\eqref{eq:sw_disp} and solving its implicit counterpart with a desired $\theta$, we obtain the dispersion relationship in $\theta\in[-\pi/2,\pi/2]$. The comparison between the numerical (circle marks) and analytical (dashed line) results are shown in Fig.~\ref{fig:2D_4}. The results match each other well.

\subsection{3D periodicity} \label{ssec:3d}

\begin{figure}[htbp] \label{fig:3D}
    \subfigure[]{
    \includegraphics[width=0.30\textwidth]{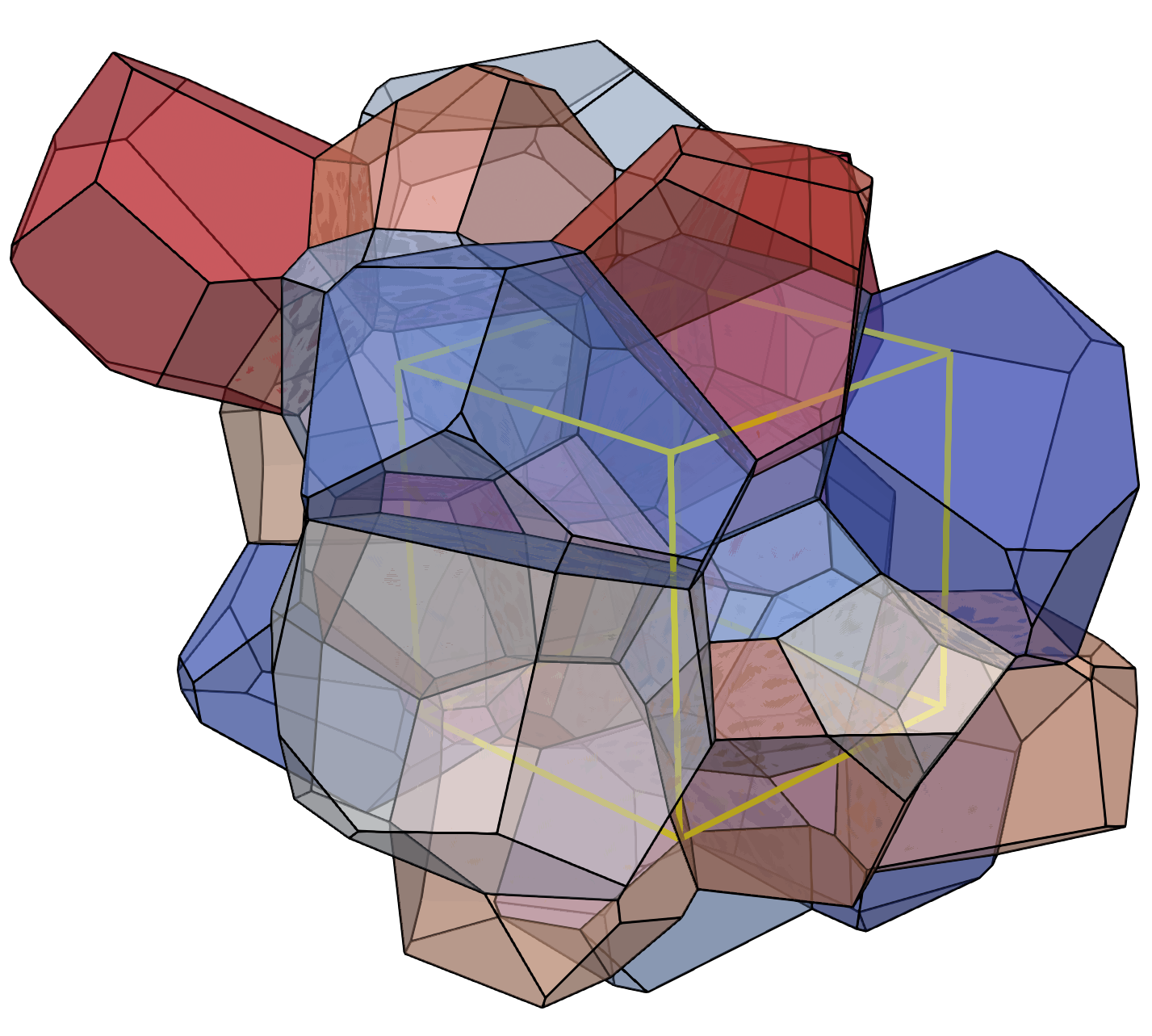} \label{fig:3D_1}
    }
    \subfigure[]{
    \includegraphics[width=0.30\textwidth]{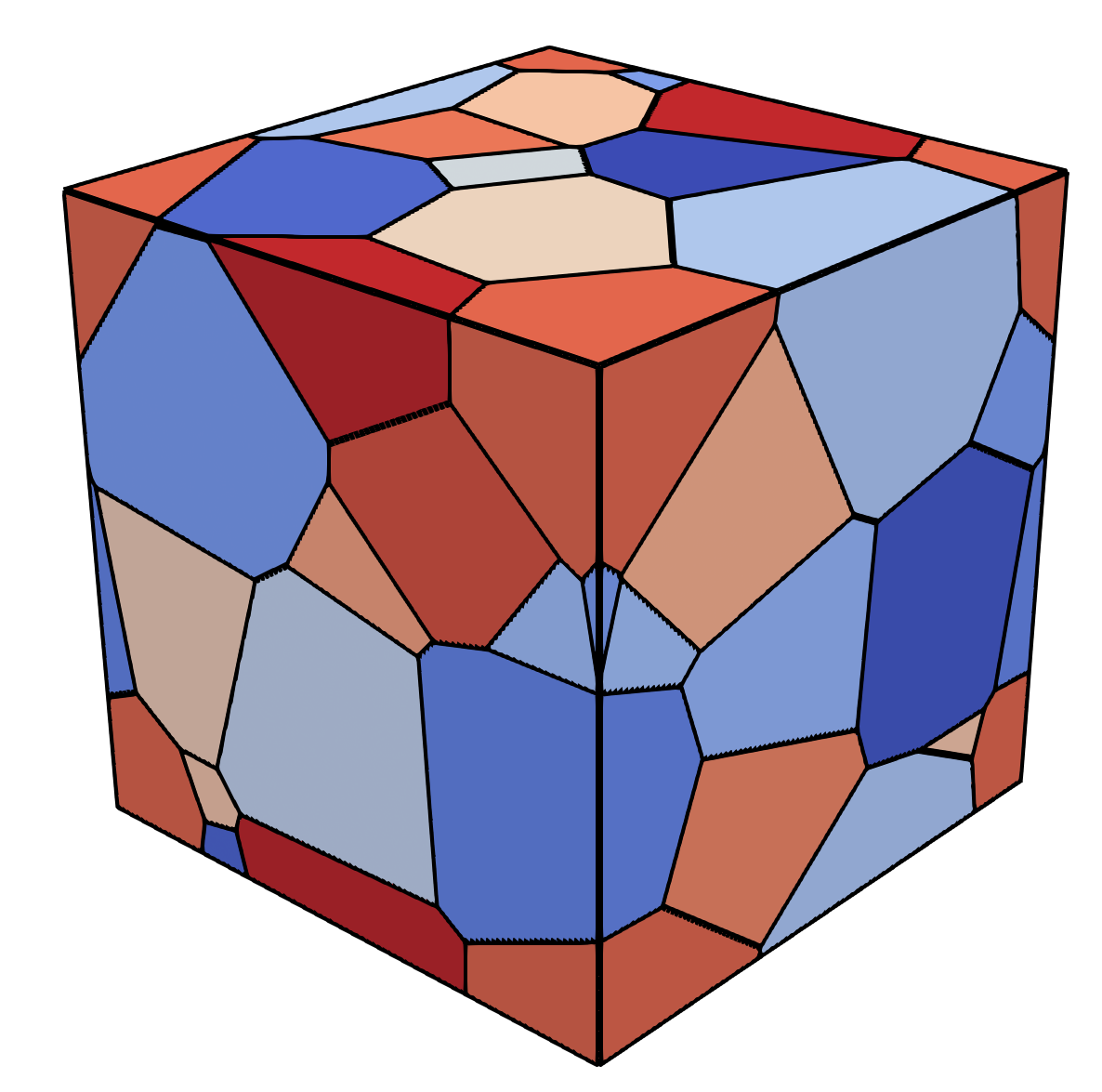} \label{fig:3D_2}
    }
    \subfigure[]{
    \includegraphics[width=0.36\linewidth]{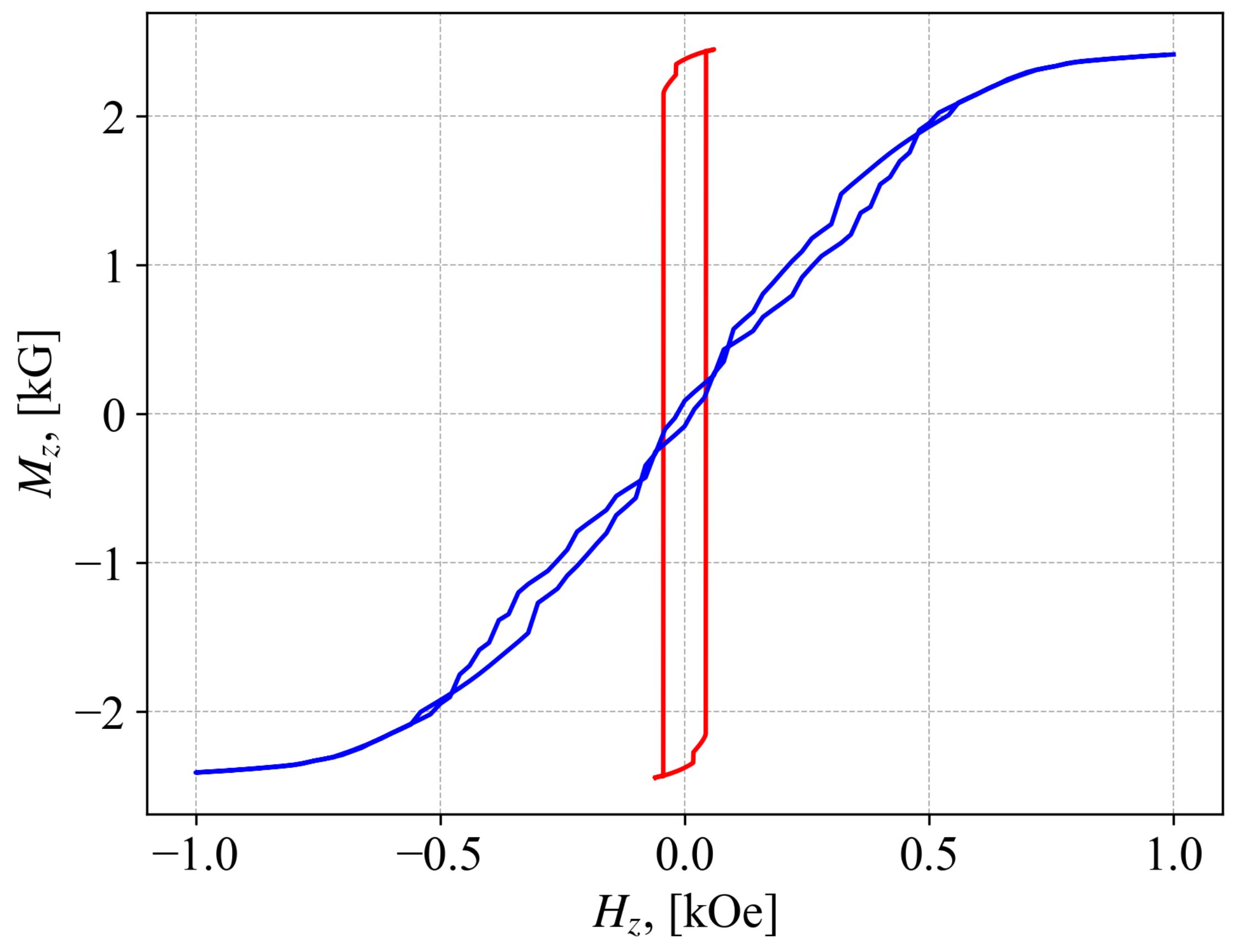} \label{fig:3D_3}
    }
    \caption{Example of 3D T-NP-PBC case.  (a)  Multi-grain structure (the unit cell is marked as a yellow cube); (b) multi-grain structure folded into the unit cell; (c) hysteresis loop for the periodic multi-grain structure. The unit cell size is $6\mu$m, the average grain size is $2.5\mu$m. The structure was meshed in a tetrahedral mesh with $70$ million elements and $12$ million nodes.}
    \label{fig:3d}
\end{figure}

For a 3D periodicity problem example, we use a complex multi-grain structure representing a granular soft magnetic material. The structure consists of 30 unique grains (Fig.~\ref{fig:3D_2}). The grains are configured to be periodic in terms of the corresponding dimensions. The structure is classified as both touching and protruding (Fig.~\ref{fig:3D_1}). After folding the protduding parts, the resulting unit cell is a cube with 3D PBC. The side length of the unit cell is $Dx = Dy = Dz = 6~\mu m$. The unit cell has $M_s=200~\mathrm{emu/cm^3}$,  $A_{ex}=1\times 10^{-6}~\mathrm{erg/cm}$, and cubic anisotropy of $K_1=3\times 10^{4}~\mathrm{erg/cm^3}$ with a random axis direction per grain. Additionally, the grains are separated by a small (nm-length) randomly variable distance and they have an interfacial exchange coupling with the surface energy density of $A_{iex}=0.05~\mathrm{erg/cm^2}$. The corresponding exchange length is $\lambda_{ex}=\sqrt{A_{ex}/2\pi M_s^2}\approx 20~\text{nm}$, and the mesh size is set to 0.5 $\lambda_{ex}$.  We impose the 3D T-NP-PBC along all directions and set $L_{[x,y,z]}$ to be same as $D_{[x,y,z]}$ to mimic an infinite domain and eliminate undesired magnetostatic field edge effects that are present without PBC. The simulation is to obtain a $M-H$ hysteresis loop. As seen from Fig.~\ref{fig:3D_3}, the non-PBC case has a significantly lower remanence magnetization and higher saturation field due to the magnetic charges at the boundary of the unit cell, which signifies the necessity of applying the PBC condition in soft material analysis.

\subsection{Performance} \label{ssec:perf}

\begin{figure}
    \centering
    \includegraphics[width=0.6\linewidth]{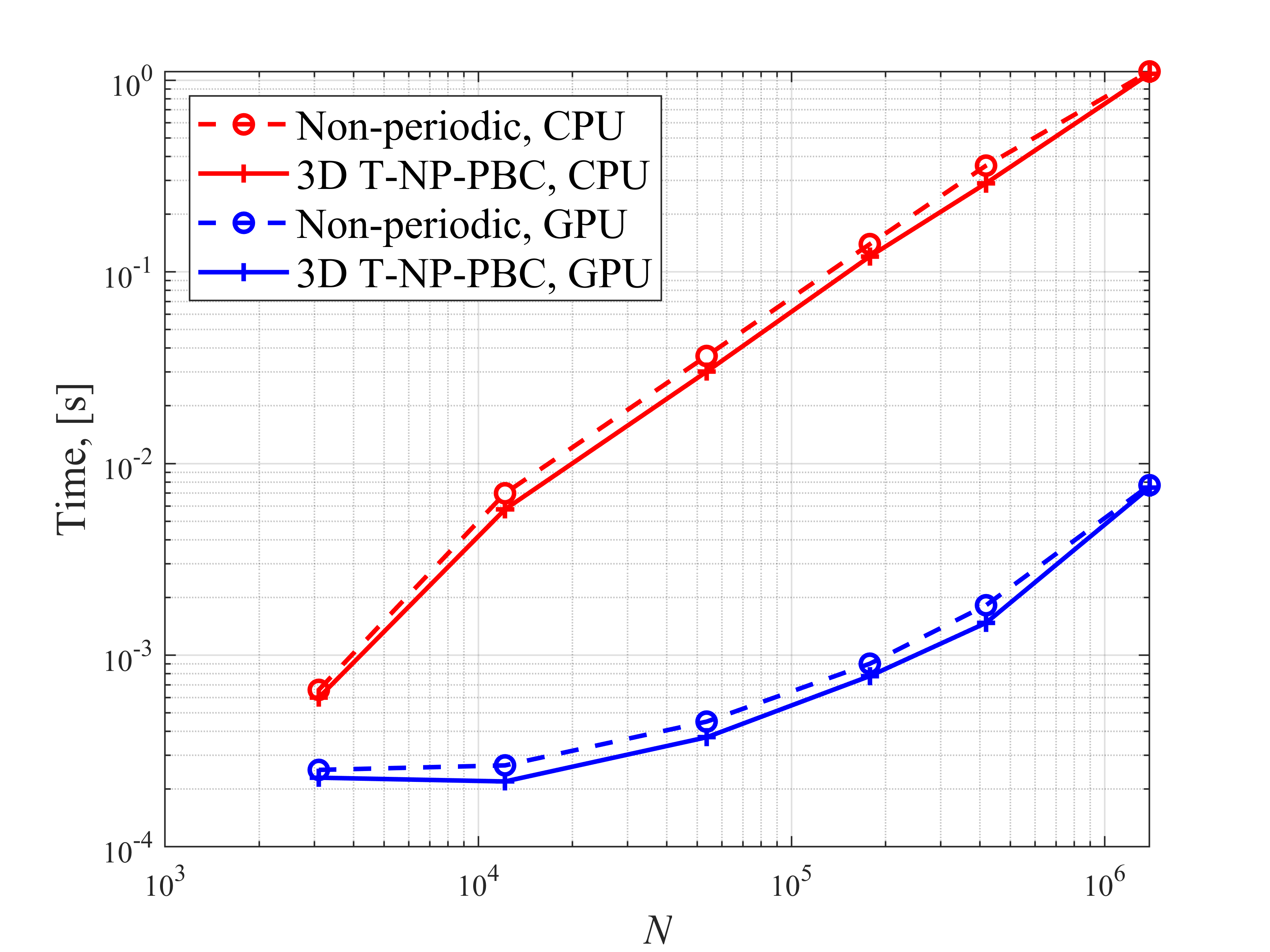}
    \caption{Performance comparison of calculation of magnetostatic field between the PM-FEM on 3D T-NP-PBC case and non-periodic case with original method on single-thread CPU and GPU.}
    \label{fig:perf}
\end{figure}

To demonstrate the computational performance, Fig.~\ref{fig:perf} shows the computational time of a one-time evaluation of the  magnetostatic field for the PM-FEM and its corresponding non-periodic case with the original BAIM-based method. The magnetostatic field evaluation is shown since it is often one of the important computational bottlenecks. To best demonstrate the performance, we choose to compare the 3D T-NP-PBC case that possesses maximal complexity as 1D/2D cases require a lower computational time. The unit cell is a cube of size $D_x=D_y=D_z=100$ with periodicity settings of $L_x=D_x$, $L_y=D_y$, $L_z=D_z$. We mesh the cube with a tetrahedral mesh. By changing the tetrahedron size, we generate models with different numbers of vertices $N$. We use the PM-FEM method and the original non-PBC method on one core of AMD 5950X CPU and on NVIDIA RTX 3080 Ti GPU with single precision and relative error level of $10^{-3}$. The results show that the periodic case may be even faster than the non-periodic case. The performance gain is from the step 2 in which the FFT-based convolution is utilized to calculate the potential on the grid. Thanks to PGF, we avoid the extra zero-padding and mirroring, which leads to a higher performance. The computational time on CPU scales nearly linearly with $N$. On GPU, for smaller $N$, the computational time scales less than linearly because the GPUs are underutilized. For larger sizes, with a full GPU utilization, the computational time scaling is also linear with $N$.  

\begin{figure}
    \centering
    \includegraphics[width=0.55\linewidth]{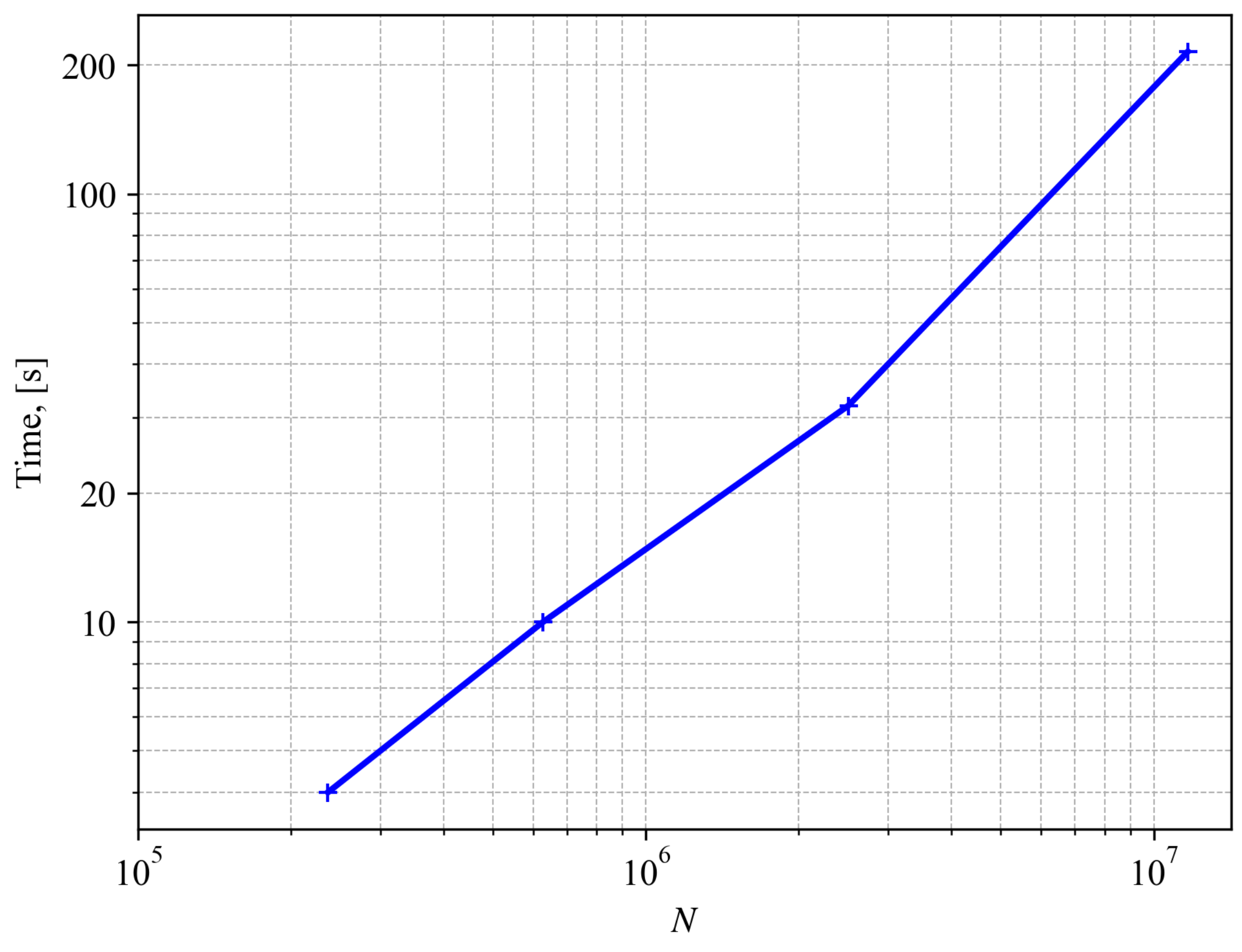}
    \caption{Performance of 1 ns simulation on single GPU.}
    \label{fig:perf_llg}
\end{figure}

Finally, Fig.~\ref{fig:perf_llg} shows the computational time of computing the magnetization dynamics in the 3D PBC case of the granular structure of Fig.~\ref{fig:3D_1}, which can be used to calculate the core loss and permeability of soft magnetic materials. The excitation was by a uniformly applied magnetic field of $1.0$ Oe at the frequency of $5$ MHz. The computational time is given for a $1\text{ns}$ simulation. The simulations were done on a NVIDIA A100 GPU. For the $1~\text{ns}$ simulation, it took $101$ time steps, $112$ non-linear iterations, and $888$ linear iterations used for the implicit BDF time stepping. The results show a mostly linear computational time scaling with $N$. For the largest considered case of $N=11.7\;\mathrm{M}$, the computational time is $216\;s$, which allows doing a design of experiment for such complex structures in real time. The computational time for a non-periodic case for the same structure was nearly the same. The achievable limit on $N$ is set by the available GPU memory, which is related not only to the periodic code components, such as the exchange and magnetostatic field, but also to other components, such as preconditioners for the linear solver used for implicit time stepping as well as data structure required for outputs.    

\section{Summary} \label{sec:con}

We introduced an efficient PM-FEM approach to account for periodicity in micromagnetic simulations. PM-FEM is applicable to all 1D/2D/3D periodic problems within a uniform framework. PM-FEM implementation is based on the non-periodic FastMag micromagnetic framework, which is updated in several aspects to include PBC. The main modifications are in handling the exchange and magnetostatic fields. Handling the exchange field is based on modifying the construction of the sparse matrix representing the Laplace operator to include elements corresponding the periodic extension in the case of touching periodic unit cells. Handling the magnetostatic field includes the local operators, such as the gradient, divergence, and surface magnetic charges as well as the long range superposition operator. The local operators for the magnetostatic field computation are handled similar to the Laplace operator for the exchange field. The long-range superposition operator involved a rapidly convergent approach for computing the PGF as well as an updated FFT-based BAIM approach to compute the PSP in $O(N\log N)$ computational cost. The modifications allows for the BAIM extension are based on the identification of the local error correction range that is based on the periodic extensions of the sources from the proper sides of the periodic unit cell.

Numerical examples demonstrate the efficiency and generality of PM-FEM. It is shown that FM-FEM can handle any 1D, 2D, and 3D periodicities with the speed that is comparable or even higher than that for similar non-periodic problems. The results also show that the updated FastMag simulator can handle large complex meshes of tens of millions of elements. PM-FEM can be used for many micromagnetic problems, including those requiring large meshes. Examples of using PM-FEM with FastMag are the study of the magnetization dynamics in periodic structures, mimicking infinite domains, such as wires, films, and bulks, and the study of spin wave propagation. 

\section{Data availability}
Data will be made available on request.

\section{Acknowledgments}
This work was supported in part by the Quantum Materials for Energy Efficient Neuromorphic-Computing (Q-MEEN-C), an Energy Frontier Research Center funded by the U.S. Department of Energy, Office of Science, Basic Energy Sciences under Award No. DESC0019273. The work was also supported in part by Binational Science Foundation, grant \#2022346. The work used Purdue Anvil cluster at Rosen Center for Advanced Computing (RCAC) in Purdue University and Delta cluster at National Center for
Supercomputing Applications (NCSA) in University of Illinois Urbana-Champaign through allocation ASC200042 from the Advanced Cyberinfrastructure Coordination Ecosystem: Services \& Support (ACCESS) program \cite{10.1145/3569951.3597559}, which is supported by National Science Foundation grants \#2138259, \#2138286, \#2138307, \#2137603, and \#2138296.

\bibliographystyle{ieeetr}
\bibliography{references}

\end{document}